\documentclass[letterpaper, 10pt, conference]{ieeeconf}
\pdfoutput=1\IEEEoverridecommandlockouts 
\usepackage{amsmath,amssymb,color,url}
\usepackage{graphicx,subfigure,warmread}
\usepackage[all,import]{xy}

\newcommand{\norm}[1]{\ensuremath{\left\| #1 \right\|}}

\newcommand{\bracket}[1]{\ensuremath{\left[ #1 \right]}}

\newcommand{\refeqn}[1]{(\ref{eqn:#1})}
\newcommand{\reffig}[1]{Fig. \ref{fig:#1}}
\newcommand{\tr}[1]{\mbox{tr}\ensuremath{\negthickspace\bracket{#1}}}
\newcommand{\deriv}[2]{\ensuremath{\frac{\partial #1}{\partial #2}}}
\newcommand{\G}{\ensuremath{\mathsf{G}}}
\newcommand{\SO}{\ensuremath{\mathsf{SO(3)}}}
\newcommand{\T}{\ensuremath{\mathsf{T}}}
\renewcommand{\L}{\ensuremath{\mathsf{L}}}
\newcommand{\so}{\ensuremath{\mathfrak{so}(3)}}
\newcommand{\SE}{\ensuremath{\mathsf{SE(3)}}}

\renewcommand{\Re}{\ensuremath{\mathbb{R}}}

\newcommand{\D}{\ensuremath{\mathbf{D}}}

\newcommand{\Ad}{\ensuremath{\mathrm{Ad}}}

\newcommand{\g}{\ensuremath{\mathfrak{g}}}

\renewcommand{\sb}{\ensuremath{\overline{s}}}
\newcommand{\mub}{\ensuremath{\overline{\mu}}}

\newcommand{\EditML}[1]{{\color{red}\protect #1}}
\renewcommand{\EditML}[1]{{\protect #1}}

\makeatletter
\def\url@leostyle{%
  \@ifundefined{selectfont}{\def\UrlFont{\sf}}{\def\UrlFont{\small\ttfamily}}}
\makeatother
\title{\LARGE \bf
Geometric Numerical Integration for\\ Complex Dynamics of Tethered Spacecraft}

\author{Taeyoung Lee\authorrefmark{1}, Melvin Leok\authorrefmark{2}, and N. Harris McClamroch%
\thanks{Taeyoung Lee, Mechanical and Aerospace Engineering, Florida Institute of Technology, Melbourne, FL 39201 {\tt taeyoung@fit.edu}}%
\thanks{Melvin Leok, Mathematics, University of California at San Diego, La Jolla, CA 92093 {\tt mleok@math.ucsd.edu}}%
\thanks{N. Harris McClamroch, Aerospace Engineering, University of Michigan, Ann Arbor, MI 48109 {\tt
nhm@umich.edu}}%
\thanks{\textsuperscript{\footnotesize\ensuremath{*}}This research has been supported in part by NSF under grants CMMI-1029551.}
\thanks{\textsuperscript{\footnotesize\ensuremath{\dagger}}This research has been supported in part by NSF under grants DMS-0714223, DMS-0726263, DMS-0747659.}
}

\begin{document}
\allowdisplaybreaks
\maketitle \thispagestyle{empty} \pagestyle{empty}

\begin{abstract}
This paper presents an analytical model and a geometric numerical integrator for a tethered spacecraft model that is composed of two rigid bodies connected by an elastic tether. This model includes important dynamic characteristics of tethered spacecraft in orbit, namely the nonlinear coupling between deformations of tether, rotational dynamics of rigid bodies, a reeling mechanism, and orbital dynamics. A geometric numerical integrator, referred to as a Lie group variational integrator, is developed to numerically preserve the Hamiltonian structure of the presented model and its Lie group configuration manifold. The structure-preserving properties are particularly useful for studying complex dynamics of a tethered spacecraft over a long period of time. These properties are illustrated by numerical simulations.
\end{abstract}

\section{Introduction}

Tethered spacecraft are composed of multiple satellites in orbit, that are connected by a thin, long cable. Numerous innovative space missions have been envisaged, such as propulsion by momentum exchange, extracting energy from the Earth's magnetic field, satellite de-orbiting, or Mars exploration~\cite{CosLor97,LorGroAA90,LesBruAA02}, and several actual missions, such as TSS, SEDS, or YES2 by NASA and ESA~\cite{CosLor97,empty0}.

The dynamics of tethered spacecraft involves nonlinear coupling effects between several dynamic modes evolving on multiple length and time scales. For example, the length of tether typically varies from $20\,\mathrm{km}$ to $100\,\mathrm{km}$, but the orbital radius of tethered spacecraft is several thousand kilometers. The natural frequency of the tether is much higher compared to the rotational attitude dynamics or the orbital period of the spacecraft. The rotational dynamics of spacecraft is nontrivially coupled to the tension of the tether, which is affected by the reeling mechanism and orbital maneuver. Therefore, it is important to accurately model tether dynamics, attitude dynamics of spacecraft, reeling mechanisms, gravitational force and the interaction between them.

Several analytic and numerical models have been developed for tethered spacecraft. However, due to the complexities of tethered spacecraft, it is common practice to use simplified models. Two point masses connected by a rigid tether is considered in~\cite{WilJSR06}. A massless, flexible tether dynamics is included in~\cite{BelLev93}. Transverse vibrations of two point masses connected by a flexible, but inextensible, tether are studied in~\cite{SomIesJGCD05}. These simplified models allow for rigorous mathematical analysis, but they may fail to predict the behaviors of an actual tethered spacecraft accurately, particularly given the fact that tethered spacecraft operations are based on weak nonlinear effects over a long time period. Recent numerical studies consider more sophisticated tethered spacecraft models including a varying tether length. But, in these advanced models, rigid body dynamics is ignored~\cite{SteZemAA95,KruPotND06}, and a reeling mechanism is neglected~\cite{SteZemAA95,ManAgrTA05}.

The goal of this paper is to develop a high-fidelity analytical model and numerical simulations for tethered spacecraft. This is an extension of preliminary work that studies a string pendulum model with a reeling mechanism~\cite{LeeLeoND09,LeeLeoPICDC09}. The first part of this paper provides a realistic and accurate analytical tethered spacecraft model including tether deformations, attitude dynamics of rigid bodies, and a reeling mechanism. We show that the governing equations of motion can be developed using Hamilton's principle. 

The second part of this paper deals with a geometric numerical integrator for tethered spacecraft. Geometric numerical integration is concerned with developing numerical integrators that preserve geometric features of a system, such as invariants, symmetry, and reversibility~\cite{HaiLub00}. 
A geometric numerical integrator, referred to as a Lie group variational integrator, has been developed for a Hamiltonian system on an arbitrary Lie group in~\cite{Lee08}.

A tethered spacecraft is a Hamiltonian system, and its configuration manifold is expressed as the product of the Lie groups $\SO$, $\SE$, and the space of connected curve segments on $\Re^3$. This paper develops a Lie group variational integrator for tethered spacecraft based on the results presented in~\cite{Lee08}. The proposed geometric numerical integrator preserves symplecticity and momentum maps, and exhibits desirable energy conservation properties. It also respects the Lie group structure of the configuration manifold, and avoids the singularities and computational complexities associated with the use of local coordinates. It can be used to study non-local, large amplitude and deformation maneuvers of tethered spacecraft accurately over a long time period. 



\section{Tethered Spacecraft}\label{sec:TS}

\begin{figure}
\centerline{\subfigure[Reference configuration]{
\renewcommand{\xyWARMinclude}[1]{\includegraphics[height=0.65\columnwidth]{#1}}
{\footnotesize\selectfont
$$\begin{xy}
\xyWARMprocessEPS{RefConf}{pdf}
\xyMarkedImport{}
\xyMarkedMathPoints{1-15}
\end{xy}
$$}}
\subfigure[Deformed configuration]{
{\footnotesize\selectfont
\renewcommand{\xyWARMinclude}[1]{\includegraphics[height=0.65\columnwidth]{#1}}
$$\begin{xy}
\xyWARMprocessEPS{DefConf}{pdf}
\xyMarkedImport{}
\xyMarkedMathPoints{1-15}
\end{xy}
$$}}}
\caption{Tethered spacecraft model}\label{fig:TSM}\vspace*{-0.4cm}
\end{figure}
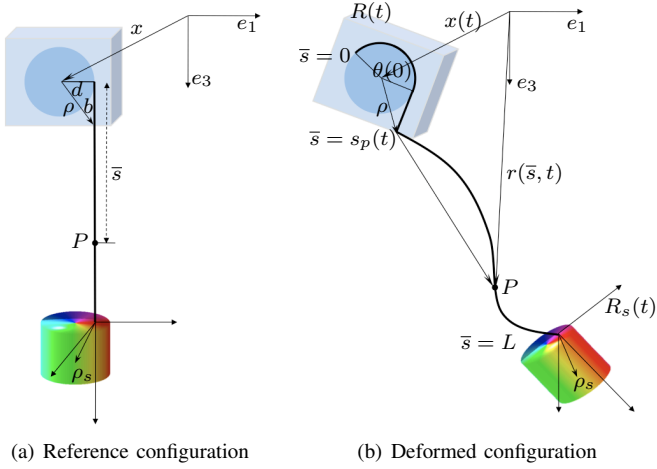

We consider two rigid spacecraft connected by an elastic tether. We assume that rigid spacecraft can freely translate and rotate in a three-dimensional space, and the tether is extensible and flexible. The bending stiffness of the tether is not considered as the diameter of the tether is assumed to be negligible compared to its length. The tether is connected to a reeling drum in a base spacecraft, and the other end of the tether is connected to a sub-spacecraft. The point where the tether is attached to the spacecraft is displaced from the center of mass so that the dynamics of the spacecraft is coupled to the tether deformations and displacements. This model is illustrated in \reffig{TSM}.

We choose a global reference frame and two body-fixed frames. The global reference frame is located at the center of the Earth. The first body fixed frame is located at the center of mass of the base spacecraft, and the second  the body-fixed frame is located at the end of the tether where the tether is attached to the sub-spacecraft. Since the tether is extensible, we need to distinguish between the arc length for the stretched deformed configuration and the arc length for the unstretched reference configuration. Define

\begin{center}
\begin{tabular}{lp{6.0cm}}
$m\in\Re$ & the mass of the base spacecraft \\
$J\in\Re^{3\times 3}$ & the inertia matrix of the base spacecraft \\
$R\in\SO$ & the rotation matrix from the first body fixed frame to the reference frame \\
$\Omega\in\Re^3$ & the angular velocity of the base spacecraft represented in its body fixed frame \\
$x\in\Re^3$ & the location of the center of mass of the base spacecraft represented in the global reference frame \\
$d\in\Re$ & the radius of the reeling drum \\
$b\in\Re$ & the length of the guideway \\
$\rho\in\Re^3$ & the vector from the center of mass of the base spacecraft to the beginning of the guideway represented in its body fixed frame, $\rho=[d,0,b]$. \\
$m_r\in\Re$ & the mass of the reeling drum \\
$J_r=\in\Re^{3\times 3}$ & the inertia matrix of the reeling drum, $J_r=\kappa_r d^2$ for a matrix $\kappa_r\in\Re^{3\times 3}$\\
$L\in\Re$ & the total unstretched length of the tether \\
$\overline s\in[0,L]$ & the unstretched arc length of the tether between the point at which the tether is attached to the reeling drum an a material point $P$ on the tether \\
\end{tabular}
\begin{tabular}{lp{6.0cm}}
$s(\overline s,t)\in\Re^+$ & the stretched arc length of the tether to the material point located at $\overline s$ \\
$s_p(t)\in[b,L]$ & the arc length of the tether between the point at which the tether is attached to the reeling drum and the beginning of the guide way \\
$r(\overline s,t)\in\Re^3$ & the deformed location of a material point $P$ from the origin of the global reference frame; $r(s_p,t)=x(t)+ R(t)\rho$ \\
$\theta(\sb)\in\Re$ & $\theta=(s_p-b-\sb)/d$ for $\sb\in[0,s_p-b]$ \\
$\overline\mu\in\Re$ & The mass of the tether per unit unstretched length \\
$m_s\in\Re$ & the mass of the sub-spacecraft \\
$J_s\in\Re$ & the inertia matrix of the sub-spacecraft \\
$R_s\in\SO$ & the rotation matrix from the second body fixed frame to the global reference frame \\
$\Omega_s\in\Re^3$ & the angular velocity of the sub-spacecraft  represented in its body fixed frame \\
$\rho_s\in\Re^3$ & the vector from the point where the tether is attached to the sub-spacecraft to the center of mass of the sub-spacecraft represented in its body fixed frame\\
$u\in\Re$ & control moment applied at the reeling drum\\
\end{tabular}
\end{center}

A configuration of this system can be described by the locations of all the material points of the tether, $r(\sb,t)$ for $\sb\in[0,L]$, the location of the base spacecraft, the attitude of both spacecraft, and the length of the deployed portion of the tether. So, the configuration manifold is $\G=C^\infty([0,l],\Re^3)\times\SE\times\SO\times\Re$, where $C^\infty([0,l],\Re^3)$ denotes the space of smooth connected curve segments on $\Re^3$, $\SO=\{R\in\Re^{3\times 3}\,|\, R^TR=I, \det[R]=1\}$, and $\SE=\Re^3\textcircled{s}\SO$~\cite{MarRat99}.

Throughout this paper, we assume that: (i) the radius of a reeling drum and the length of a guideway is small compared to the length of a tether; (ii) the reeling drum rotates about the second axis of the first body fixed frame attached to the base spacecraft; (iii) the deployed portion of the tether is extensible, but the portion of the tether on the reel and the guideway inside of the base spacecraft is inextensible; (iv) the gravity is uniform over the base spacecraft and the sub-spacecraft.

\section{Continuous-time Analytical Model}\label{sec:AM}

In this section, we develop continuous-time equations of motion for a tethered spacecraft using Hamilton's variational principle. The attitude kinematics equation of the base spacecraft and the sub-spacecraft is given by
\begin{align}
    \dot R = R\hat\Omega,\quad \dot R_s = R\hat\Omega_s,\label{eqn:Rdot}
\end{align}
where the \textit{hat map} $\hat\cdot:\Re^3\rightarrow\so$ is defined by the condition that $\hat x y=x\times y$ for any $x,y\in\Re^3$. 

\subsection{Lagrangian}

\paragraph*{Kinetic energy}

The kinetic energy of the base spacecraft excluding the reeling drum is given by
\begin{align}
T_{b_1}=\frac{1}{2}m \dot x \cdot \dot x + \frac{1}{2} \Omega \cdot J\Omega.\label{eqn:Tb1}
\end{align}
Under the assumption that the radius of a reeling drum is much less than the length of the tether, the kinetic energy of the reeling drum and the part of the tether inside of the base spacecraft can be approximated by
\begin{align}
T_{b_2}= \frac{1}{2}(m_r+\mub s_p) \dot x \cdot \dot x + \frac{1}{2}\mub s_p \dot s_p^2
+\frac{1}{2} \kappa_2 \dot s_p^2.\label{eqn:Tb2}
\end{align}
where $\kappa_2 = e_2\cdot \kappa_r e_2$. Let $\dot r(\sb,t)$ be the partial derivative of $r(\sb,t)$ with respect to $t$. The kinetic energy of the deployed portion of the tether is given by
\begin{align}
T_{t} & = \int_{s_p}^{L} \frac{1}{2} \mub \dot r(\sb)\cdot \dot r(\sb)\, d\sb.\label{eqn:Tt}
\end{align}
Let $\tilde\rho\in\Re^3$ be the vector from the end of the tether to a mass element of the sub-spacecraft represented with respect to its body fixed frame. The location of the mass element in the global reference  frame is given by $r(L)+R_s\tilde\rho$. Then, the kinetic energy of the sub-spacecraft is given by
\begin{align}
T_{s} & = \int_{\mathcal{B}_s} \frac{1}{2} \|\dot r(L)+R_s\hat\Omega_s\tilde\rho        \|^2\,dm\nonumber\\
&=\frac{1}{2} m_s \dot r(L)\cdot\dot r(L) + m_s\dot r(L)\cdot R_s\hat\Omega_s \rho_s + \frac{1}{2}\Omega_s \cdot J_s\Omega_s.\label{eqn:Ts}
\end{align}
Here, we use the fact that $\int_{\mathcal{B}_s} \tilde \rho dm = \rho_c$, and $J_s = -\int_{\mathcal{B}_s}{\hat\rho_s}^2 dm$. The total kinetic energy is given by $T=T_{b_1}+T_{b_2}+T_t+T_s$. 

\paragraph*{Potential energy}
By the assumption that the size of reeling drum is small compared to the length of the tether, the gravitational potential of the base spacecraft and the reeling mechanism is approximated as follows
\begin{align}
V_{b} 
& = -(m+m_r+\mub s_p)\frac{GM}{\|x\|},\label{eqn:Vb}
\end{align}
where the gravitational constant and the mass of the Earth are denoted by $G$ and $M$, respectively. 

The strain of the tether at a material point located at $r(\sb)$ is given by
\begin{align*}
\epsilon = \lim_{\Delta \sb\rightarrow 0} \frac{\Delta s(\sb)-\Delta \sb}{\Delta \sb}
= s'(\sb) -1,
\end{align*}
where $(\;)'$ denote the partial derivative with respect to $\sb$. The tangent vector at the material point is given by 
\begin{align*}
e_t= \deriv{r(\sb)}{s} = \deriv{r(\sb)}{\sb}\deriv{\sb}{s(\sb)} = \frac{r'(\sb)}{s'(\sb)}.
\end{align*}
Since this tangent vector has unit length, we have $s'(\sb)=\norm{r'(\sb)}$. Therefore, the strain can be written as $\epsilon = \norm{r'(\sb)}-1$. Using this, the elastic potential and the gravitational potential of the deployed portion of the tether is given by
\begin{align}
V_t & = 
\frac{1}{2}\int_{s_p}^L EA(\norm{r'(\sb)}-1)^2\,d\sb
-\int_{s_p}^L \mub  \frac{GM}{\|r(\sb)\|}\,d\sb,\label{eqn:Vt}
\end{align}
where $E$ and $A$ denote the Young's modulus and the sectional area of the tether. 

The location of the center of mass of the sub-spacecraft is $r(L)+R_s\rho_s$ in the global reference frame. Since we assume that the gravity is uniform over each spacecraft body, the gravitational potential energy of the sub-spacecraft is given by
\begin{align}
V_s = -m_s \frac{GM}{\|r(L)+R_s\rho_s\|}.\label{eqn:Vs}
\end{align}

From \refeqn{Tb1}-\refeqn{Vs}, the Lagrangian of the tethered spacecraft is given by
\begin{align}
L = T_{b_1} + T_{b_2} + T_t + T_s - V_b - V_t - V_s.\label{eqn:L}
\end{align}

\subsection{Euler-Lagrange Equations}
Let the action integral be $\mathfrak{G} = \int_{t_0}^{t_f} L\,dt$. According to Hamilton's principle, the variation of the action integral is equal to the negative of the virtual work for fixed boundary conditions, which yields Euler-Lagrange equations. For the given tethered spacecraft model, this requires the following three careful consideration: (i) the domain of the integral depends on the variable $s_p(t)$ at \refeqn{Tt}; (ii) the rotation matrices $R,R_s$ that represents the attitudes lie in the nonlinear Lie group $\SO$; (iii) as the tether is assumed to be inextensible in the guideway, and it is extensible outsize of the guideway, there exists a discontinuity in strain at the beginning of the guideway. 

\paragraph{Time-Varying Domain} 

Due to \refeqn{Tt}, the variation of the action integral $\delta\mathfrak{G}$ includes the following term, $\int_{t_0}^{t_t}\int_{s_p}^L \mub \dot r(\sb)\cdot \delta\dot r(\sb)\, d\sb dt$. Here, we cannot apply integration by parts with respect to $\sb$, since the order of the integrals cannot be interchanged due to the time dependency in the variable $s_p(t)$.  Instead, we use Green's theorem~\cite{LeeLeoND09},
\begin{align}
\oint_\mathcal{C} \dot r(\sb)\cdot\delta r(\sb) \,d\sb 
= \int_{t_0}^{t_f} \int_{s_p(t)}^L \frac{d}{dt}(\dot r(\sb)\cdot\delta r(\sb))\,d\sb dt,\label{eqn:Green1}
\end{align}
where $\oint_\mathcal{C}$ represents the counterclockwise line integral on the boundary $\mathcal{C}$ of the region $[t_0,t_f]\times[s_p(t),L]$. The boundary $\mathcal{C}$ is composed of four lines: $(t=t_0,\sb\in[s_p(t_0),L])$, $(t=t_f,\sb\in[s_p(t_f),L])$, $(t\in[t_0,t_f],\sb=s_p(t))$, and $(t\in[t_0,t_f],\sb=L)$. For the first two lines, $\delta r(\sb)=0$ since $t=t_0,t_f$. For the last line, $d\sb=0$ since $\sb$ is fixed. Thus, parameterizing the third line by $t$, we obtain 
\begin{align*}
\oint_\mathcal{C} \dot r(\sb)\cdot\delta r(\sb) \,d\sb 
=\int_{t_0}^{t_f} \dot r(s_p(t))\cdot \delta r(s_p(t))\, \dot s_p(t)\, dt.
\end{align*}
Substituting this into \refeqn{Green1} and rearranging, we obtain
\begin{align}
&\int_{t_0}^{t_f}  \int_{s_p}^L \dot r(\sb)\cdot \delta \dot r(\sb)\,d\sb dt\nonumber\\
& =\int_{t_0}^{t_f} \bracket{\int_{s_p}^L -\ddot r(\sb)\cdot \delta r(\sb) \,d\sb 
+ \dot r(s_p)\cdot \delta r(s_p)\, \dot s_p\,} dt.\label{eqn:drdr}
\end{align}

\begin{figure*}[!t]
\small
\setcounter{equation}{15}
\begin{gather}
-(m+m_r+\mub s_p)\ddot x  -GM(m+m_r+\mub s_p) \frac{x}{\|x\|^3}+\mub\dot s_p(-r'(s_p^+)\dot s_p - R\hat\rho\Omega)+F(s_p)=0,\label{eqn:EL0}\\
-J\dot\Omega -\hat\Omega J\Omega+\mub\dot s_p\hat\rho R^T(-r'(s_p^+)\dot s_p+\dot x - R\hat\rho\Omega)+\hat\rho R^T F(s_p)-ue_2=0,\label{eqn:EL1}\\
-(\mub s_p + \kappa_2)\ddot s_p-\frac{1}{2}\mub (\dot x -R\hat\rho\Omega)\cdot (\dot x -R\hat\rho\Omega)+\frac{1}{2}\mub \dot x \cdot \dot x
-\mub \frac{GM}{\| r(s_p)\|}+\mub \frac{GM}{\|x\|}
-F(s_p)\cdot r'(s_p^+)+\frac{u}{d}=0,\label{eqn:EL3}\\
-\mub\ddot r(\sb)+F'(\sb)-\mub GM \frac{r(\sb)}{\|r(\sb)\|}=0,\qquad (\sb\in[s_p,L],\quad r(s_p)=x+R\rho),\label{eqn:EL4}\\
-m_s\ddot r(L)+m_sR_s\hat\rho_s \dot\Omega_s-m_sR_s\hat\Omega_s^2\rho_s-m_sR_s\hat{\dot\Omega}_s\rho_s -GMm_s \frac{r(L)+R_s\rho_s}{\|r(L)+R_s\rho_s\|^3}-F(L)=0,\label{eqn:EL5}\\
-J_s\dot\Omega_s-m_s\hat\rho_s R_s^T\ddot r(L)
-\hat\Omega_s J_s\Omega_s-GMm_s \hat\rho_s R_s^T\frac{r(L)+R_s\rho_s}{\|r(L)+R_s\rho_s\|^3}=0,\label{eqn:EL6}
\end{gather}
where $F(\sb)=EA\frac{\norm{r'(\sb)}-1}{\norm{r'(\sb)}}r'(\sb)$ denotes the tension of the tether.

\setcounter{equation}{11}
\hrulefill
\end{figure*}

\paragraph{Variation of Rotation Matrices}

The attitudes of spacecraft are represented by the rotation matrix $R,R_s\in\SO$. Therefore, the variation of the rotation matrix should be consistent with the geometry of the special orthogonal group. In~\cite{Lee08}, it is expressed in terms of the exponential map as
\begin{align}
\delta R =\frac{d}{d\epsilon}\bigg|_{\epsilon=0} R^\epsilon= \frac{d}{d\epsilon}\bigg|_{\epsilon=0} R \exp \epsilon\hat\eta = R \hat\eta, \label{eqn:delR}
\end{align}
for $\eta\in\Re^3$. The key idea is expressing the variation of a Lie group element in terms of a Lie algebra element. This is desirable since the Lie algebra $\so$ of the special orthogonal group, represented by $3\times 3$ skew-symmetric matrices, is isomorphic as a Lie algebra to $\Re^3$. As a result, the variation of the three-dimenstional rotation matrix $R$ is expressed in terms of a vector $\eta\in\Re^3$. We can directly show that \refeqn{delR} satisfies $\delta(R^TR)=\delta R^T R + R^T \delta R =-\hat\eta+\hat\eta=0$. The corresponding variation of the angular velocity is obtained from the kinematics equation \refeqn{Rdot}:
\begin{align}
\delta\hat\Omega = \frac{d}{d\epsilon}\bigg|_{\epsilon=0} (R^\epsilon)^T \dot R^\epsilon = (\dot\eta + \Omega\times \eta)^\wedge.\label{eqn:delw}
\end{align}

\paragraph{Variational Principle with Discontinuity}

Let $r(s_p^-)$, and $r(s_p^+)$ be the material point of the tether just inside the guide way, and the material point just outside the guide way, respectively. Since the tether is inextensible inside the guide way, $\|r'(s_p^-)\|=1$. Since the tether is extensible outside the guide way, $\|r'(s_p^+)\|=1+\epsilon^+$, where $\epsilon^+$ represents the strain of the tether just outside the guide way. Due to this discontinuity, the speed of the tether changes instantaneously by the amount $\epsilon^+|\dot s_p|$ at the guide way.

As a result, the variation of the action integral is not equal to the negative of the virtual work done by the external control moment $u$ at the reeling drum. Instead, an additional term $Q$, referred to as Carnot energy loss term should be introduced \cite{KruPotND06,CreJanJAM97}. The resulting variational principle is given by
\begin{align}
\delta\mathfrak G + \int_{t_0}^{t_f} ( Q + u/d )\delta s_p-ue_2\cdot\eta \, dt =0.\label{eqn:VP}
\end{align}
The corresponding time rate of change of the total energy is given by $\dot E = (Q + u/d ) \dot s_p$, where the first term $Q\dot s_p$ represents the energy dissipation rate due to the velocity and strain discontinuity. The corresponding expression for $Q$ has been developed in~\cite{LeeLeoND09}:
\begin{align}
Q = -\frac{1}{2}\mub (\|r'(s_p^+)\|-1)^2 \dot s_p^2 - \frac{1}{2} EA (\|r'(s_p^+)\|-1)^2 .
\end{align}

\paragraph{Euler-Lagrange Equations} Using these results, and the variational principle with discontinuity \refeqn{VP}, we obtain the Euler-Lagrange equations for the given tethered spacecraft model in \refeqn{EL0}-\refeqn{EL6}. In \refeqn{EL4}, we require that $r(s_p)=x+R\rho$ for the continuity of the tether.

These can be simplified in a number of special cases. For example, we can substitute $\dot s_p=0$ when the length of the deployed portion of the tether is fixed, and we can set $\rho=\rho_s=0$ when the main spacecraft and the sub-spacecraft are modeled as point masses instead of rigid bodies. 

\setcounter{equation}{21}

\section{Lie Group Variational Integrator}\label{sec:LGVI}

The Euler-Lagrange equations developed in the previous section provide an analytical model for a tethered spacecraft. However, the standard finite difference approximations or finite element approximations of those equations using a general purpose numerical integrator may not preserve the geometric properties of the system accurately~\cite{HaiLub00}.

Lie group variational integrators provide a systematic method of developing geometric numerical integrators for Lagrangian/Hamiltonian systems evolving on a Lie group~\cite{Lee08}. As they are  derived from a discrete analogue of Hamilton's principle, they preserve symplecticity and the momentum map, and it exhibits good total energy behavior. They also preserve the Lie group structure as they update a group element using the group operation. These properties are critical for accurate and efficient simulations of complex dynamics of multibody systems~\cite{LeeLeoCMDA07}. 

In this section, we develop a Lie group variational integrator for a tethered spacecraft.

\subsection{Discretized Tethered Spacecraft  Model}

Let $h>0$ be a fixed time-step. The value of variables at $t=t_0+kh$ is denoted by a subscript $k$.  We discretize the deployed portion of the tether using $N$ identical line elements. Since the unstretched length of the deployed portion of the tether is $L-s_{p_k}$, the unstretched length of each element is $l_{k} = \frac{L-s_{p_k}}{N}$. Let the subscript $a$ denote the variables related to the $a$-th element. The natural coordinate on the $a$-th element is defined by
\begin{align}
\zeta_{k,a} (\sb) = \frac{(\sb-s_{p_k})-(a-1)l_{k}}{l_{k}}
\end{align}
for $\sb\in[s_{p_k}+(a-1)l_{k},s_{p_k}+al_{k}]$. This varies between 0 and 1 on the $a$-th element. Let $S_0,S_1$ be shape functions given by $S_0(\zeta)=1-\zeta$, and $S_1(\zeta)=\zeta$. These shape functions are also referred to as  \textit{tent functions}.

Using this finite element model, the position vector $r(\sb,t)$ of a material point in the $a$-th element is approximated as follows:
\begin{align}
r_k(\sb)= S_0(\zeta_{k,a}) r_{k,a} + S_1(\zeta_{k,a}) r_{k,a+1}.\label{eqn:rksb}
\end{align}

Therefore, a configuration of the presented discretized tethered spacecraft at $t=kh+t_0$ is described by $g_k=(x_k;R_k;s_{p_k};r_{k,1},\ldots,r_{k,N+1};R_{s_k})$, and the corresponding configuration manifold is $\G=\Re^3\times\SO\times\Re\times (\Re^3)^{N+1}\times \SO$. This is a Lie group where the group acts on itself by the diagonal action~\cite{MarRat99}: the group action on $x_k$, $s_{p_k}$, and $r_{k,a}$ is addition, and the group action on $R_k,R_{s_k}$ is matrix multiplication.

We define a discrete-time kinematics equation using the group action as follows. Define $f_k=(\Delta x_k;F_k;\Delta s_{p_k}; \Delta r_{k,1},\ldots$, $\Delta r_{k,N+1};F_{s_k})\in\G$ such that $g_{k+1}=g_k f_k$:
\begin{align}
    &(x_{k+1};R_{k+1};s_{p_{k+1}};r_{k+1,a};R_{k+1})=\nonumber\\
    &(x_k+\Delta x_k;R_kF_k;s_{p_k}+\Delta s_{p_k};r_{k,a}+\Delta r_{k,a};R_{s_k} F_{s_k}).
\end{align}
Therefore, $f_k\in \G$ represents the relative update between two integration steps. This ensures that the structure of the Lie group configuration manifold is numerically preserved since $g_{k}$ is updated by $f_k$ using the right Lie group action of $\G$ on itself.

\subsection{Discrete Lagrangian}

A discrete Lagrangian $L_d(g_k,f_k):\G\times\G\rightarrow\Re$ is an approximation of the Jacobi solution of the Hamilton--Jacobi equation, which is given by the integral of the Lagrangian along the exact solution of the Euler-Lagrange equations over a single time-step:
\begin{align*}
    L_d(g_k,f_k)\approx \int_0^h L(\tilde g(t),{\tilde g}^{-1}(t)\dot{\tilde g} (t))\,dt,
\end{align*}
where $\tilde g(t):[0,h]\rightarrow \G$ satisfies Euler-Lagrange equations with boundary conditions $\tilde{g}(0)=g_k$, $\tilde{g}(h)=g_kf_k$. The resulting discrete-time Lagrangian system, referred to as a variational integrator, approximates the Euler-Lagrange equations to the same order of accuracy as the discrete Lagrangian approximates the Jacobi solution~\cite{MarWesAN01}.

We construct a discrete Lagrangian for the tethered spacecraft using the trapezoidal rule.
From the attitude kinetics equations \refeqn{Rdot}, the angular velocity is approximated by
\begin{align*}
\hat\Omega_k \approx \frac{1}{h} R_k^T (R_{k+1}-R_k) = \frac{1}{h} (F_k -I).
\end{align*}
Define a non-standard inertia matrix $J_d = \frac{1}{2}\mathrm{tr}[J]I-J$. Using the trace operation and the non-standard inertia matrix, the rotational kinetic energy of the main spacecraft at \refeqn{Tb1} can be written in terms of $\hat\Omega$ as $\frac{1}{2}\Omega\cdot J\Omega=\frac{1}{2}\mathrm{tr}[\hat\Omega J_d \hat\Omega^T]$. Using this, and \refeqn{Tb1}, \refeqn{Tb2}, the kinetic energy of the base spacecraft is given by
\begin{align}
T_{k,b} & = \frac{1}{2h^2}(m+m_r+\mub s_{p_k})\Delta x_k\cdot\Delta x_k \nonumber\\
&\quad+ \frac{1}{2h^2} (\mub s_{p_k} + \kappa_s) \Delta s_{p_k}^2
+\frac{1}{h^2} \tr{(I-F_k)J_d},\label{eqn:Tkb}
\end{align}
where we use properties of the trace operator: $\mathrm{tr}[AB]=\mathrm{tr}[BA]=\mathrm{tr}[A^TB^T]$ for any matrices $A,B\in\Re^{3\times 3}$. 

Next, we find the kinetic energy of the tether. Using the chain rule, the partial derivative of $r_k(\sb)$ given by \refeqn{rksb} with respect to $t$ is given by
\begin{align*}
\dot r_k(s) & = 
\frac{1}{h} \bigg\{S_0(\zeta_{k,a}) \Delta r_{k,a} + S_1(\zeta_{k,a}) \Delta r_{k,a+1}\\
&\quad+\frac{(L-s)}{(L-s_{p_k})}\frac{(r_{k,a}-r_{k,a+1})}{l_k}\Delta s_{p_k}\bigg\}.
\end{align*}
Substituting this into \refeqn{Tt}, the contribution of the $a$-th tether element to the kinetic energy of the tether is given by
\begin{align}
T_{k,a} 
& = \frac{1}{2h^2}M^1_k \Delta r_{k,a}\cdot \Delta r_{k,a}
+\frac{1}{2h^2}M^2_k \Delta r_{k,a+1}\cdot \Delta r_{k,a+1}\nonumber\\
&+\frac{1}{2h^2}M^3_{k,a} \Delta s_{p_k}^2
+\frac{1}{h^2}M^{12}_k \Delta r_{k,a}\cdot \Delta r_{k,a+1}\nonumber\\
&+\frac{1}{h^2} M^{23}_{k,a}\Delta s_{p_k} \cdot \Delta r_{k,a+1}
+\frac{1}{h^2} M^{31}_{k,a}\Delta s_{p_k} \cdot \Delta r_{k,a},\label{eqn:Tka}
\end{align}
where inertia matrices are given by
\begin{align*}
M^1_k&=\frac{1}{3}\mub l_k,\quad 
M^2_k=M^1_k,\\
M^3_{k,a} &= \frac{1}{3}\mub l_k \frac{(3N^2+3N+1-6Na-3a+3a^2)}{N^2},
\\
M^{12}_k &= \frac{1}{6}\mub l_k,\quad
M^{23}_{k,a} = \frac{1}{6}\mub \frac{(1+3N-3a)}{N}(r_{k,a}-r_{k,a+1}),\\
M^{31}_{k,a} &= \frac{1}{6}\mub \frac{(2+3N-3a)}{N}(r_{k,a}-r_{k,a+1}).
\end{align*}
Similar to \refeqn{Tkb}, from \refeqn{Ts}, the kinetic energy of the sub-spacecraft is given by 
\begin{align}
T_{k,s} &= \frac{1}{2h^2} m_s \Delta r_{k,N+1}\cdot \Delta r_{k,N+1}+\frac{1}{h^2}\tr{(I-F_{s_k})J_{s_d}}\nonumber\\
&\quad +\frac{1}{h^2}m_s \Delta r_{k,N+1}\cdot R_{s_k}(F_{s_k}-I)\rho_s.\label{eqn:Tks}
\end{align}
From \refeqn{Tkb}, \refeqn{Tka}, \refeqn{Tks}, the total kinetic energy of the discretized tethered spacecraft is given by
\begin{align}
T_k = T_{k,b} + \sum_{a=1}^N T_{k,a} + T_{k,s}.\label{eqn:Tk}
\end{align}

Similarly, from \refeqn{Vb}, \refeqn{Vt}, and \refeqn{Vs}, the total potential energy is given by
\begin{align}
V_k & = -GM(m+m_r+\mub s_{p_k}) \frac{1}{\|x_k\|}\nonumber\\
&\quad +\sum_{a=1}^N -2 GM \mub l_k \frac{1}{\|r_{k,a}+r_{k,a+1}\|}\nonumber\\
&\quad+\frac{1}{2}\frac{EA}{l_k} (\|r_{k,a+1}-r_{k,a}\| - l_k)^2\nonumber\\
&\quad - GMm_s \frac{1}{\|r_{k,N+1}+R_{s_k}\rho_s\|}. \label{eqn:Vk}
\end{align}

Using \refeqn{Tk}, \refeqn{Vk}, we choose the discrete-Lagrangian of the discretized tethered spacecraft as follows:
\begin{align}
L_{d_k}(g_k,f_k) = h T_k(g_k,f_k) -\frac{h}{2} V_k(g_k,f_k) - \frac{h}{2} V_{k+1}(g_k,f_k).\label{eqn:Ld}
\end{align}

\subsection{Discrete-time Euler-Lagrange Equations}

We define the discrete action sum $\mathfrak{G}_d=\sum_{k=1}^n L_{d_k}(g_k,f_k)$. According to the discrete Hamilton's principle, the variation of the action sum is equal to the negative of the discrete virtual work. This yields discrete-time Euler-Lagrange equations, referred to as variational integrators. 

In~\cite{Lee08}, the following Lie group variational integrator has been developed for Lagrangian systems on an arbitrary Lie group:
\begin{gather}
\begin{aligned}
    \T_e^*\L_{f_{k-1}}\cdot & \D_{f_{k-1}}  L_{d_{k-1}}-\Ad^*_{f_{k}^{-1}}\cdot(\T_e^*\L_{f_{k}}\cdot \D_{f_{k}}L_{d_{k}})\\
    &\quad+\T_e^*\L_{g_{k}}\cdot \D_{g_{k}} L_{d_{k}}+U_{d_k}+Q_{d_k}=0,
\end{aligned}\label{eqn:DEL0}\\
    g_{k+1} = g_k f_k,\label{eqn:DEL1}
\end{gather}
where $\T\L:\T\G\rightarrow\T\G$ is the \EditML{tangent} map of the left translation, $\D_f$ represents the derivative with respect to $f$, and $\Ad^*:\G\times\g^*\rightarrow\g^*$ is $\mathrm{co}$-$\mathrm{Ad}$ operator~\cite{MarRat99}. 

The virtual work due to the control input and the Carnot energy loss are denoted by $U_{d_k}\in \g^*$, and $Q_{d_k}\in \g^*$, respectively, and they are chosen as
\begin{align}
U_{d_k}\cdot (g^{-1}\delta g_k) & = \frac{h}{d} u_{k}\delta s_k - h u_k e_2\cdot\eta_k,\label{eqn:Ud}\\
Q_{d_k}\cdot (g^{-1}\delta g_k) & = -\frac{h}{2l_k^2} (\mub \Delta s_{p_k}^2 /h^2 +EA)\nonumber\\
&\quad \times ( \norm{r_{k,2}-r_{k,1}}-l_k)^2\,\delta s_{p_k}.\label{eqn:Qd}
\end{align}

By substituting \refeqn{Ld}, \refeqn{Ud}, and \refeqn{Qd} into \refeqn{DEL0} and \refeqn{DEL1}, we obtain a Lie group variational integrator for the given discrete tethered spacecraft model. Due to page limits, we do not present the detailed results in this paper. But, we demonstrate their computational properties in the next section.

\section{Numerical Example}\label{sec:NE}

Properties of the tethered spacecraft are chosen as
\begin{gather*}
m=490\,\mathrm{kg},\quad
m_r=10\,\mathrm{kg},\quad
m_s=150\,\mathrm{kg},\\
l=120\,\mathrm{km},\quad
\mub=24.7\,\mathrm{kg/km},\quad
EA=659700\,\mathrm{N},\\
J=\mathrm{diag}[5675.8,\, 5675.8,\, 6125]\,\mathrm{kgm^2},\quad
\rho=[0.5, 0.0, 1]\mathrm{m},\\
J_s=\mathrm{diag}[500,\, 500,\, 300]\,\mathrm{kgm^2},\quad
\rho_s=[0, 0, 1]\,\mathrm{m}.
\end{gather*}
Initially, the base spacecraft is on a circular orbit with an altitude of $300\,\mathrm{km}$, and the tether and the sub-spacecraft are aligned along the radial direction. The initial unstretched length of the deployed portion of the tether is $20\,\mathrm{km}$, i.e., $s_p(0)=100\,\mathrm{km}$. The initial velocity at each point of the tether and the sub-spacecraft is chosen such that it corresponds to the velocity of a circular orbit at their altitude. 

We consider the following three cases. In the first case, the reeling drum is fixed so that the length of the deployed portion of the tether is fixed, i.e., $\dot s_p\equiv 0$. In the second case, the reeling drum is free to rotate, and the tether is released by gravity to $100\,\mathrm{km}$. The third case is the same as the first case except that the initial velocities of the base spacecraft and the sub-spacecraft are perturbed by about $15\%$ to generate a tumbling motion. These are summarized with the time-step and the simulation time as follows: 
\begin{center}\vspace*{-0.43cm}
\begin{tabular}{c|l|c|c}
& Description&  $h\,\mathrm{(s)}$ & $t_f\,\mathrm{(s)}$\\\hline
Case 1 & Fixed reeling drum & 0.05 & 6000\\
Case 2 & Releasing the tether to $100\,\mathrm{km}$ & 0.05 & 3848\\
Case 3 & Velocity perturbation of Case 1 & 0.01 & 500\\
\end{tabular}
\end{center}
Note that the orbital period of a point mass on the circular orbit with the altitude of $300\,\mathrm{km}$ is $5410$ seconds. For all cases, the number of tether elements is $N=20$.

The following figures illustrate simulation results for each case. We consider a fictitious local vertical, local horizontal (LVLH) frame that is attached to an imaginary spacecraft on a circular orbit with an altitude of $300\,\mathrm{km}$. For each figure, we have the following subfigures: (a) the maneuvers of the tethered spacecraft are illustrated with respect to the LVLH frame. To represent the attitude dynamics of spacecraft, the size of the spacecraft is increased by a factor of 100, and the relative strain distribution of the tether at each instant is represented by a color shading (animations illustrating these maneuvers are also available at \url{http://my.fit.edu/~taeyoung}). The remaining subfigures show: (b) the energy transfer, (c) the computed total energy deviation from its initial value, (d) the angular velocity of the base spacecraft, and (e) the unstretched/stretched length of the tether.

In the first case, we observe a pendulum-like motion where the tether is taut and its stretch length is almost close to the unstretched length. But, there exists a strain wave that propagates along the tether, and nontrivial attitude maneuvers for the base spacecraft and the sub-spacecraft. The proposed Lie group variational integrator exhibit excellent conservation properties: the maximum relative total energy deviation is $2.37\times10^{-8}\,\%$ of its initial value, and the maximum orthogonality error of rotation matrices is $\max\{\|I-R^TR\|\}=1.32\times 10^{-13}$.

\begin{figure}
\centerline{
	\subfigure[Snapshots observed at the LVLH frame (km) (The size of spacecraft is increased by a factor of 100 to illustrate attitude dynamics.)]{
		\includegraphics[width=0.95\columnwidth]{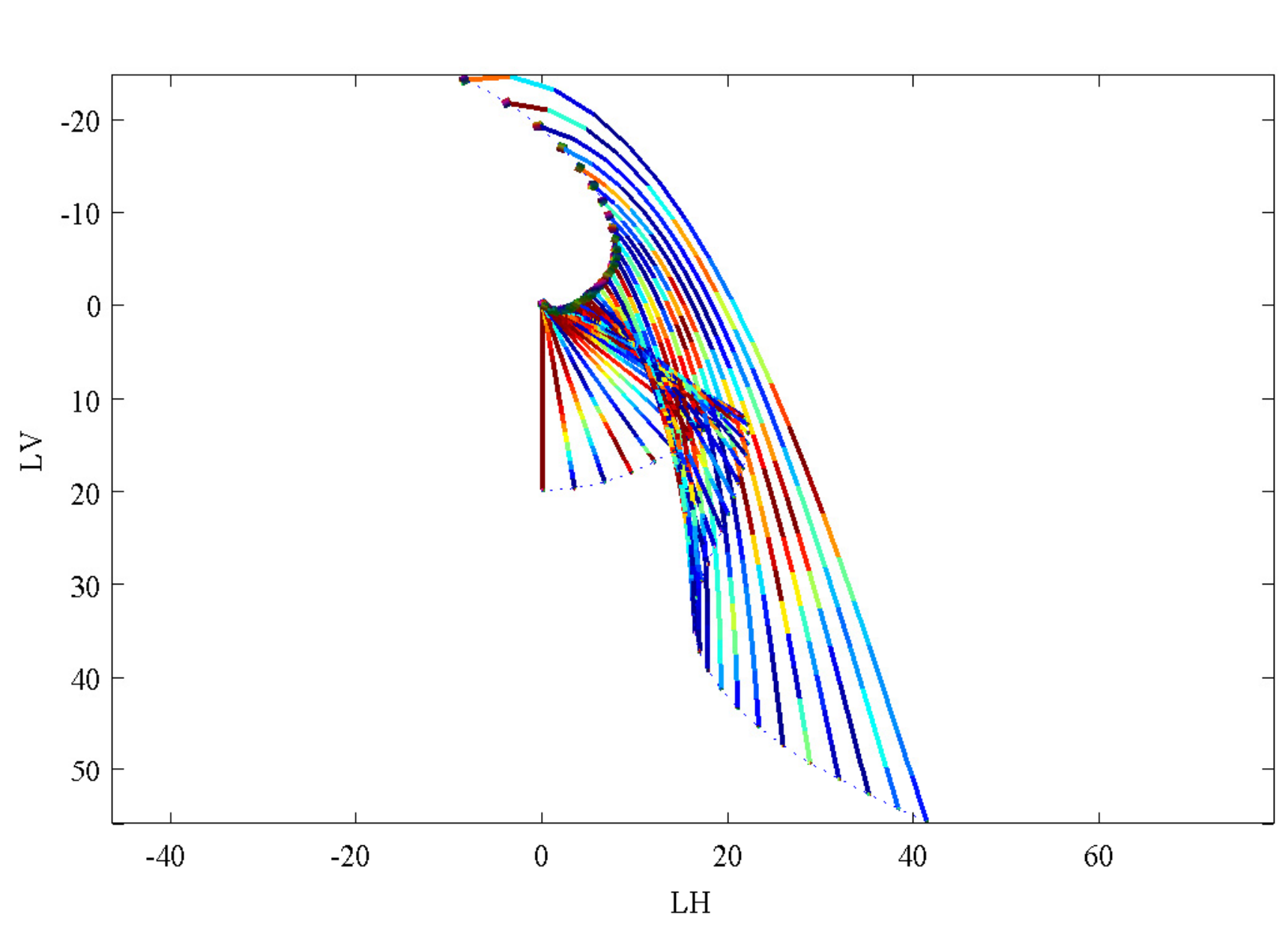}}
}
\centerline{
	\subfigure[$T_{base}+T_{sub}$ (red), $T_{tether}$ (green), $V_{gravity}$ (cyan), $V_{elastic}$ (blue), total energy (black)]{
		\includegraphics[width=0.49\columnwidth]{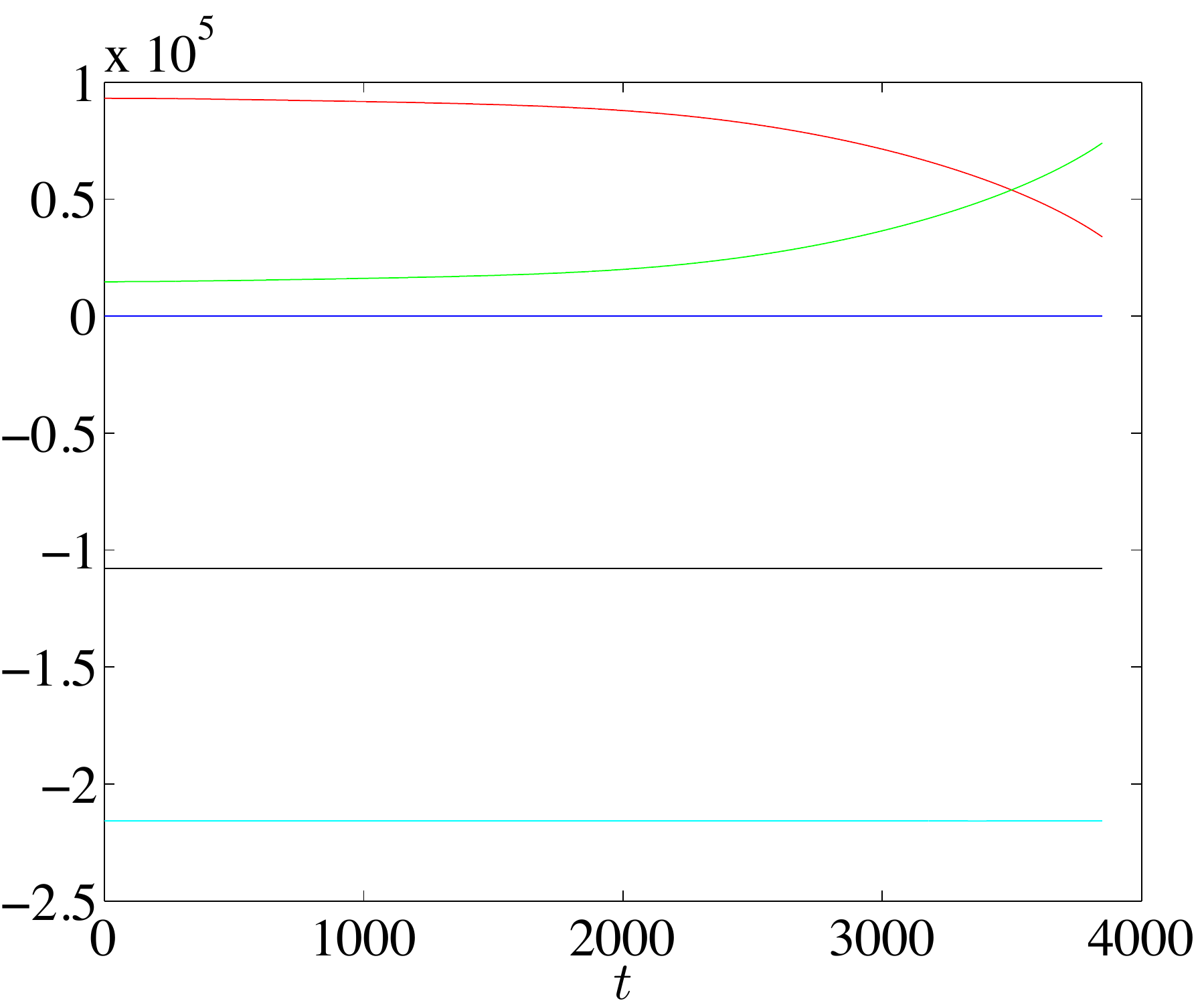}\label{fig:2Es}}
	\hspace*{0.00\textwidth}
	\subfigure[Computed total energy deviation $E(t)-E(0)$]{
		\includegraphics[width=0.51\columnwidth]{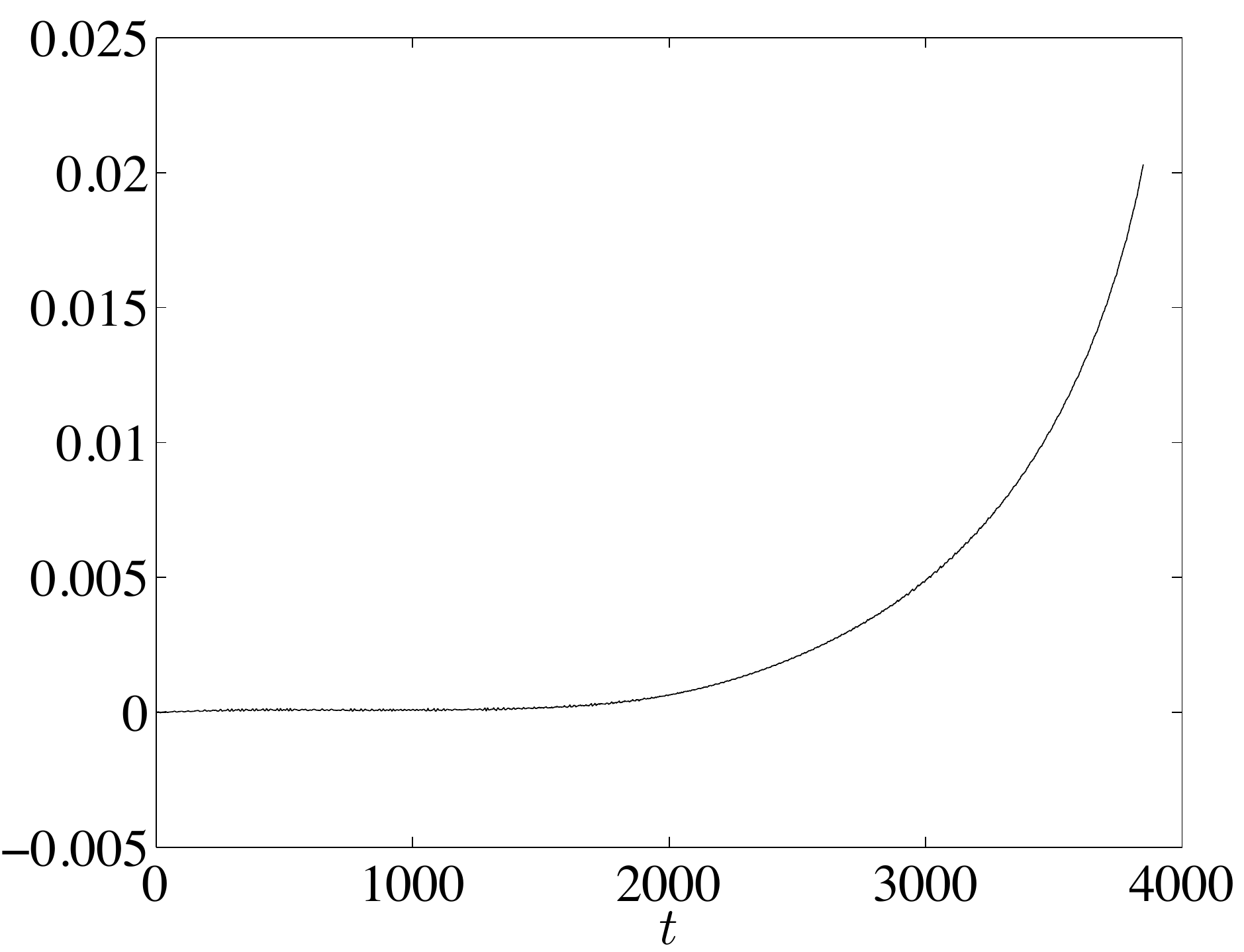}}
}
\centerline{
	\subfigure[Angular velocity of the base spacecraft $\Omega$]{
		\includegraphics[width=0.51\columnwidth]{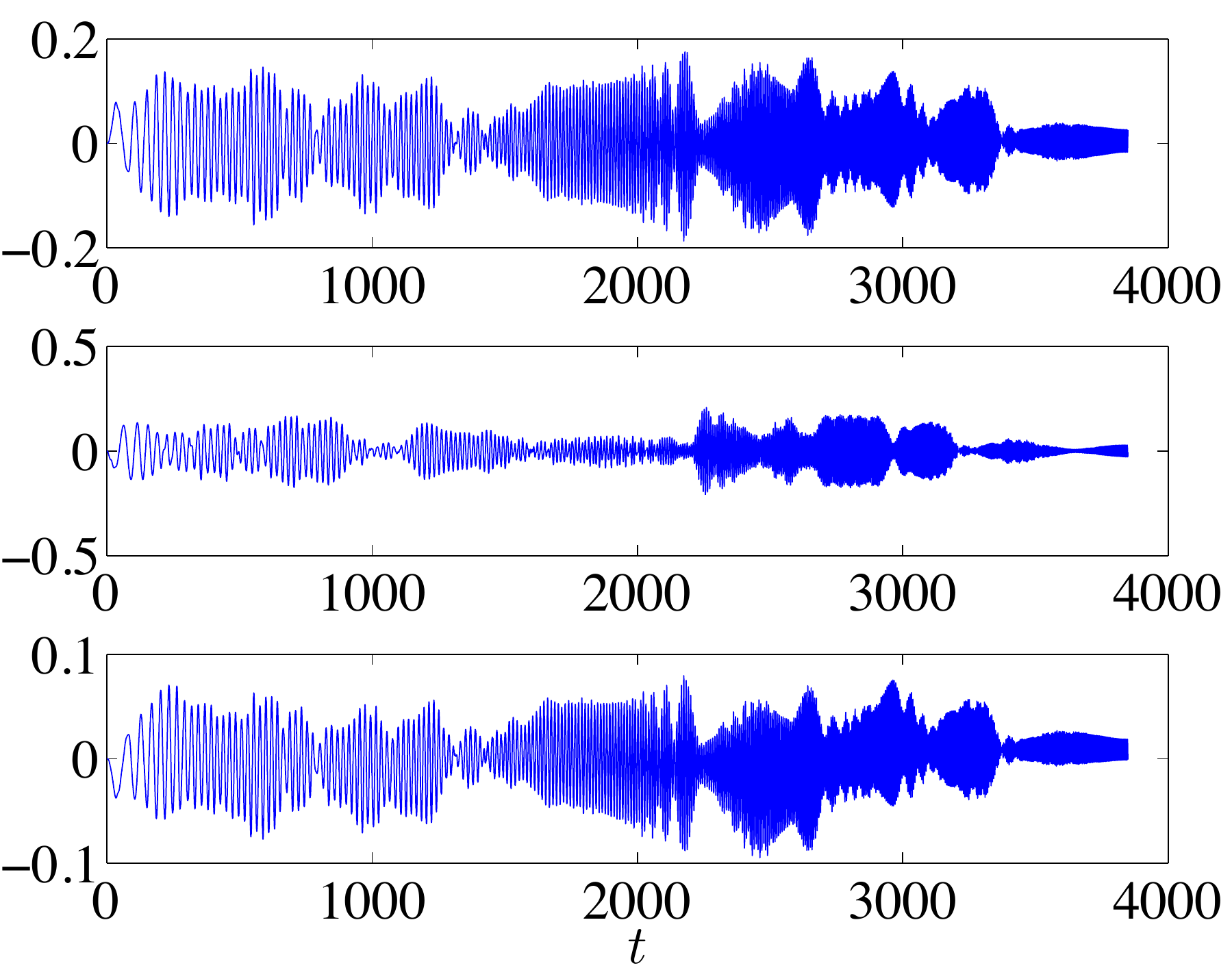}}
	\hspace*{0.00\textwidth}
	\subfigure[Unstretched length of the deployed part of the tether (red), stretched length (blue)]{
		\includegraphics[width=0.49\columnwidth]{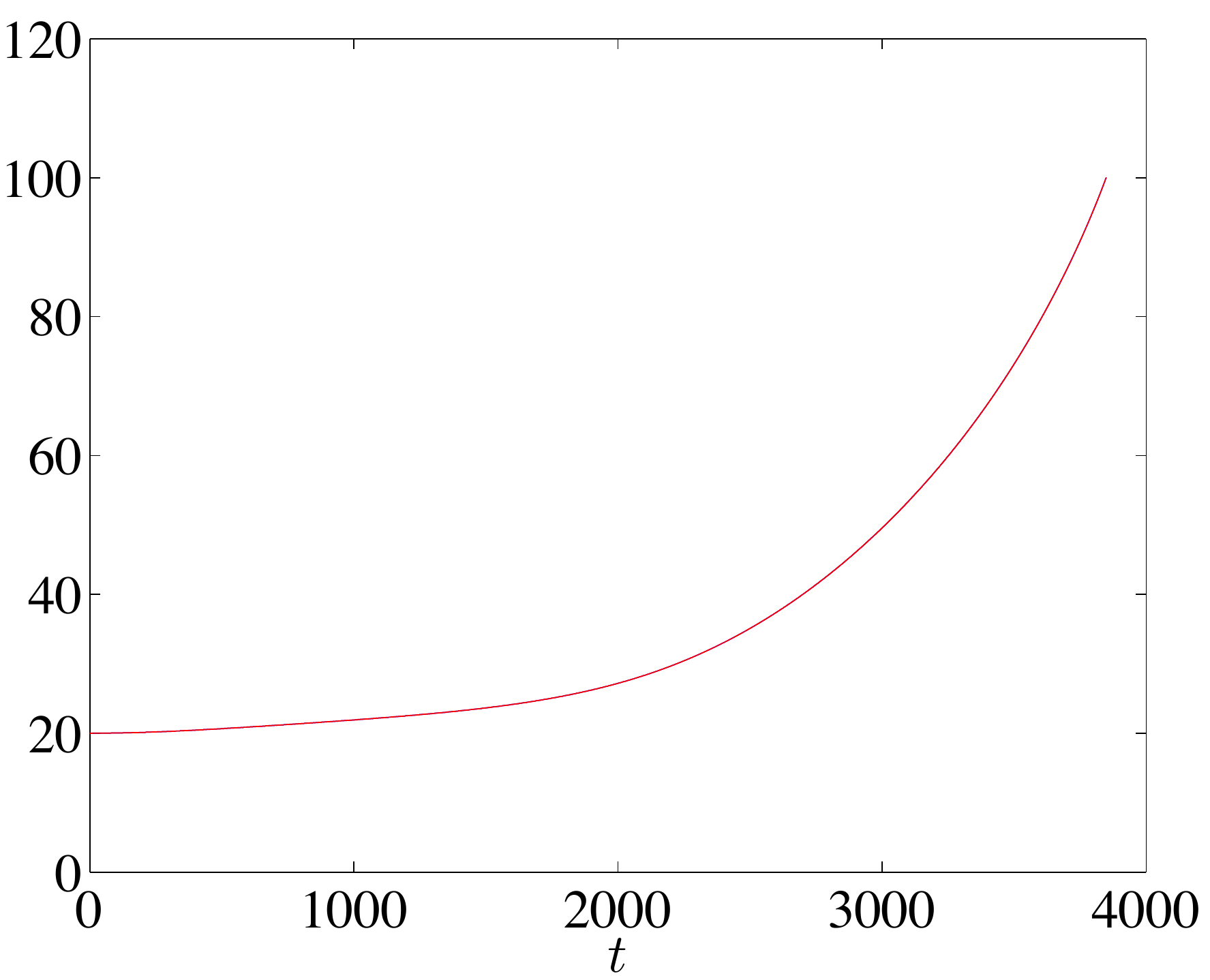}}
}

\caption{Case 2: Circular orbit, Releasing tether}
\end{figure}

In the second case, the tether is deployed by gravity gradient effects, and due to the Carnot energy term discussed in the previous section, the total energy increases slightly. As the mass in the base spacecraft is transferred to the deployed portion of the tether, there is a transfer of kinetic energy between two parts, as seen in \reffig{2Es}.

The third case is most challenging: there are in-plane and out-of-plane tumbling maneuvers, while the tether is stretched by $25\,\%$, and the attitude dynamics of spacecraft is nontrivially excited with a large angular velocity. The proposed Lie group variational integrator computes the complex dynamics of this tethered spacecraft accurately. The maximum relative total energy deviation is $3.48\times 10^{-4}\,\%$, and the maximum orthogonality error of rotation matrices is $\max\{\|I-R^TR\|\}=8.03\times 10^{-14}$.

\section{Conclusions}

We develop continuous-time equations of motion and a geometric numerical integrator for a tethered spacecraft model that includes tether deformation, spacecraft attitude dynamics, and a reeling mechanism. This provides an analytical model that is defined globally on the Lie group configuration manifold, and the Lie group variational integrator preserves the underlying geometric features, thereby yielding a reliable numerical simulation tool for complex maneuvers over a long time period.
%
%

\bibliography{ACC11}

\begin{thebibliography}{10}
\providecommand{\url}[1]{#1}
\csname url@rmstyle\endcsname
\providecommand{\newblock}{\relax}
\providecommand{\bibinfo}[2]{#2}
\providecommand\BIBentrySTDinterwordspacing{\spaceskip=0pt\relax}
\providecommand\BIBentryALTinterwordstretchfactor{4}
\providecommand\BIBentryALTinterwordspacing{\spaceskip=\fontdimen2\font plus
\BIBentryALTinterwordstretchfactor\fontdimen3\font minus
  \fontdimen4\font\relax}
\providecommand\BIBforeignlanguage[2]{{%
\expandafter\ifx\csname l@#1\endcsname\relax
\typeout{** WARNING: IEEEtran.bst: No hyphenation pattern has been}%
\typeout{** loaded for the language `#1'. Using the pattern for}%
\typeout{** the default language instead.}%
\else
\language=\csname l@#1\endcsname
\fi
#2}}

\bibitem{CosLor97}
M.~Cosmo and E.~Lorenzini, ``Tethers in space handbook,'' NASA Marshall Space
  Flight Center, Tech. Rep., 1997.

\bibitem{LorGroAA90}
E.~Lorenzini, M.~Grossi, and M.~Cosmo, ``Low altitude tethered mars probe,''
  \emph{Acta Astronautica}, vol.~21, no.~1, pp. 1--12, 1990.

\bibitem{LesBruAA02}
L.~Less, C.~Bruno, C.~Ullvieri, U.~Ponzi, M.~Parisse, G.~Laneve, G.~Vannaroni,
  M.~Dobrowolny, F.~{De Venuto}, B.~Bertotti, and L.~Anselmo, ``Satellite
  de-orbiting by means of electrodynamic tethers, {Part I}: general conecpts
  and requirements,'' \emph{Acts Astronautica}, vol.~50, no.~7, pp. 399--406,
  2002.

\bibitem{empty0}
\BIBentryALTinterwordspacing
Young engineers' satellite 2, {E}uropean {S}pace {A}gency. [Online]. Available:
  \url{http://www.esa.int/SPECIALS/YES/index.html}
\BIBentrySTDinterwordspacing

\bibitem{WilJSR06}
P.~Williams, ``Simple approach to orbital control using spinning electrodynamic
  tethers,'' \emph{Journal of Spacecraft and Rockets}, vol.~43, no.~1, pp.
  253--256, 2006.

\bibitem{BelLev93}
V.~Beletsky and E.~Levin, \emph{Dynamics of Space Tether Systems}.\hskip 1em
  plus 0.5em minus 0.4em\relax Univelt, 1993.

\bibitem{SomIesJGCD05}
L.~Somenzi, L.~Iess, and J.~Pelaez, ``Linear stability analysis of
  electrodynamic tethers,'' \emph{Journal of Guidance, Control, and Dynamics},
  vol.~28, no.~5, pp. 843--849, 2005.

\bibitem{SteZemAA95}
W.~Steiner, J.~Zemann, A.~Steindl, and H.~Troger, ``Numerical study of large
  amplitutde oscillations of a two-satellite continuous tether system with a
  varying length,'' \emph{Acta Astronautica}, vol.~35, no. 9-11, pp. 607--621,
  1995.

\bibitem{KruPotND06}
M.~Krupa, W.~Poth, M.~Schagerl, A.~Steindl, W.~Steiner, and H.~Troger,
  ``Modelling, dynamics and control of tethered satellites systems,''
  \emph{Nonlinear Dynamics}, vol.~43, pp. 73--96, 2006.

\bibitem{ManAgrTA05}
K.~Mankala and S.~Agrawal, ``Dynamic modeling and simulation of satellite
  tethered systems,'' \emph{Transactions of the ASME}, vol. 127, pp. 144--156,
  2005.

\bibitem{LeeLeoND09}
T.~Lee, M.~Leok, and N.~McClamroch, ``Computational dynamics of a 3{D} elastic
  string pendulum attached to a rigid body and an inertially fixed reel
  mechanism,'' \emph{Nonlinear Dynamics}, 2009, submitted.

\bibitem{LeeLeoPICDC09}
------, ``Dynamics of a 3{D} elastic string pendulum,'' in \emph{Proceedings of
  IEEE Conference on Decision and Control}, 2009, pp. 3347--3352.

\bibitem{HaiLub00}
E.~Hairer, C.~Lubich, and G.~Wanner, \emph{Geometric numerical integration},
  ser. Springer Series in Computational Mechanics 31.\hskip 1em plus 0.5em
  minus 0.4em\relax Springer, 2000.

\bibitem{Lee08}
T.~Lee, ``Computational geometric mechanics and control of rigid bodies,''
  Ph.D. dissertation, University of Michigan, 2008.

\bibitem{MarRat99}
J.~Marsden and T.~Ratiu, \emph{Introduction to Mechanics and Symmetry},
  2nd~ed., ser. Texts in Applied Mathematics.\hskip 1em plus 0.5em minus
  0.4em\relax Springer-Verlag, 1999, vol.~17.

\bibitem{CreJanJAM97}
E.~Crellin, F.~Janssens, D.~Poelaert, W.~Steiner, and H.~Troger, ``On balance
  and variational formulations of the equations of motion of a body deploying
  along a cable,'' \emph{Journal of Applied Mechanics}, vol.~64, pp. 369--374,
  1997.

\bibitem{LeeLeoCMDA07}
T.~Lee, M.~Leok, and N.~H. McClamroch, ``Lie group variational integrators for
  the full body problem in orbital mechanics,'' \emph{Celestial Mechanics and
  Dynamical Astronomy}, vol.~98, no.~2, pp. 121--144, June 2007.

\bibitem{MarWesAN01}
J.~Marsden and M.~West, ``Discrete mechanics and variational integrators,'' in
  \emph{Acta Numerica}.\hskip 1em plus 0.5em minus 0.4em\relax Cambridge
  University Press, 2001, vol.~10, pp. 317--514.

\end{thebibliography}
\bibliographystyle{IEEEtran}

\newpage
\begin{figure}
\centerline{
	\subfigure[Snapshots observed at the LVLH frame (km) (The size of spacecraft is increased by a factor of 100 to illustrate attitude dynamics.)]{
		\includegraphics[width=0.95\columnwidth]{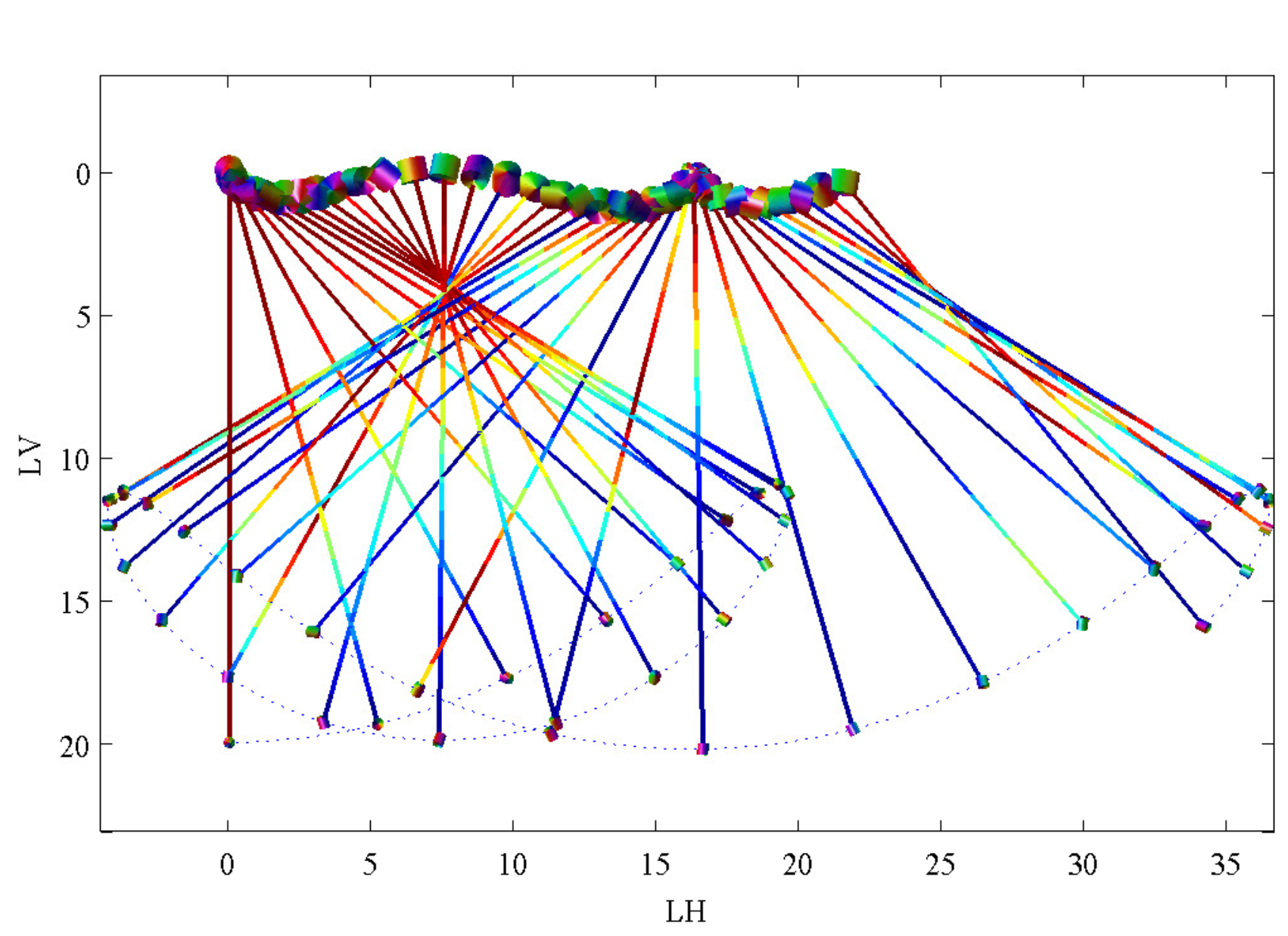}}
}
\centerline{
	\subfigure[$T_{base}+T_{sub}$ (red), $T_{tether}$ (green), $V_{gravity}$ (cyan), $V_{elastic}$ (blue), total energy (black)]{
		\includegraphics[width=0.5\columnwidth]{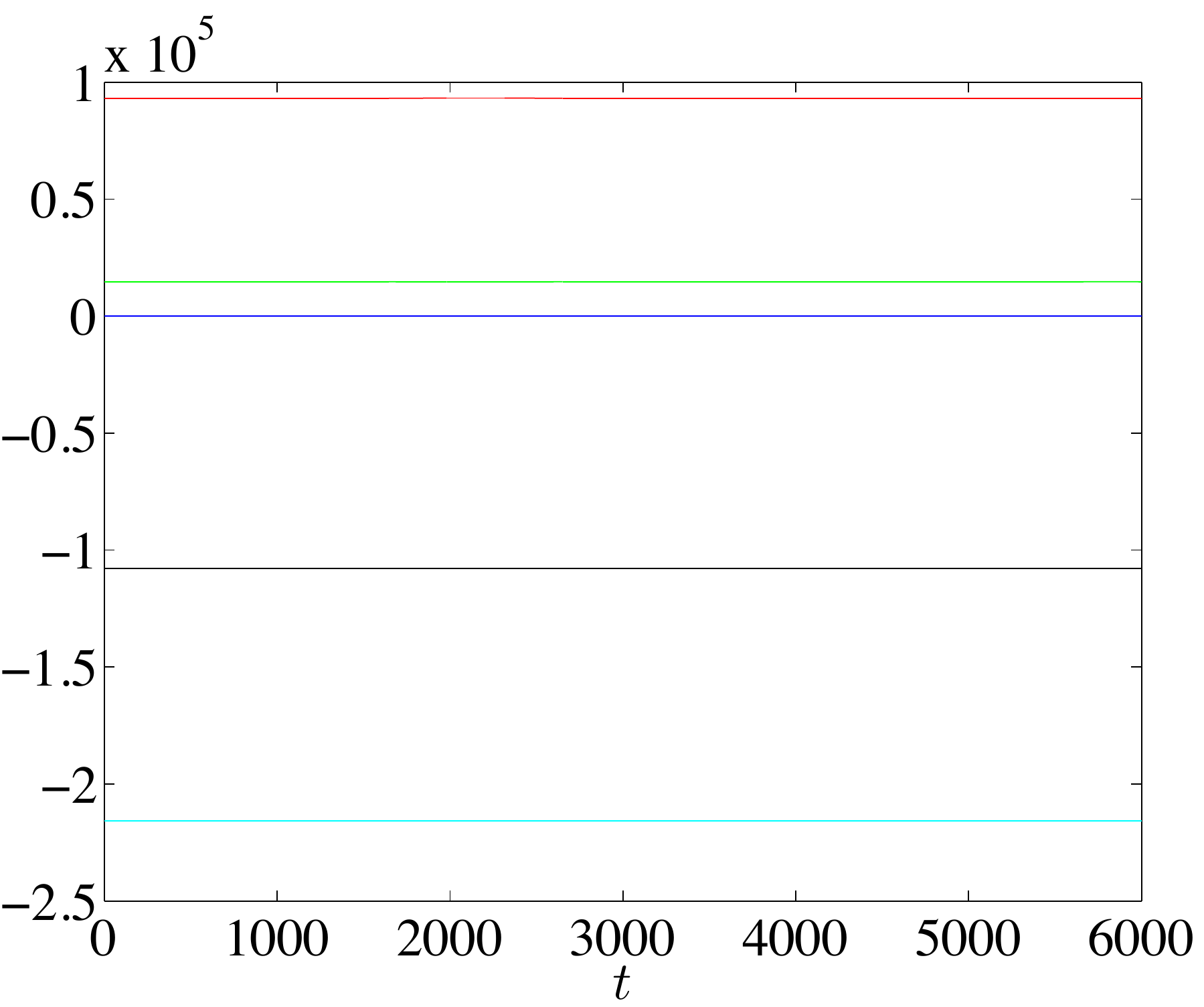}}
	\hspace*{0.00\textwidth}
	\subfigure[Computed total energy deviation $E(t)-E(0)$]{
		\includegraphics[width=0.49\columnwidth]{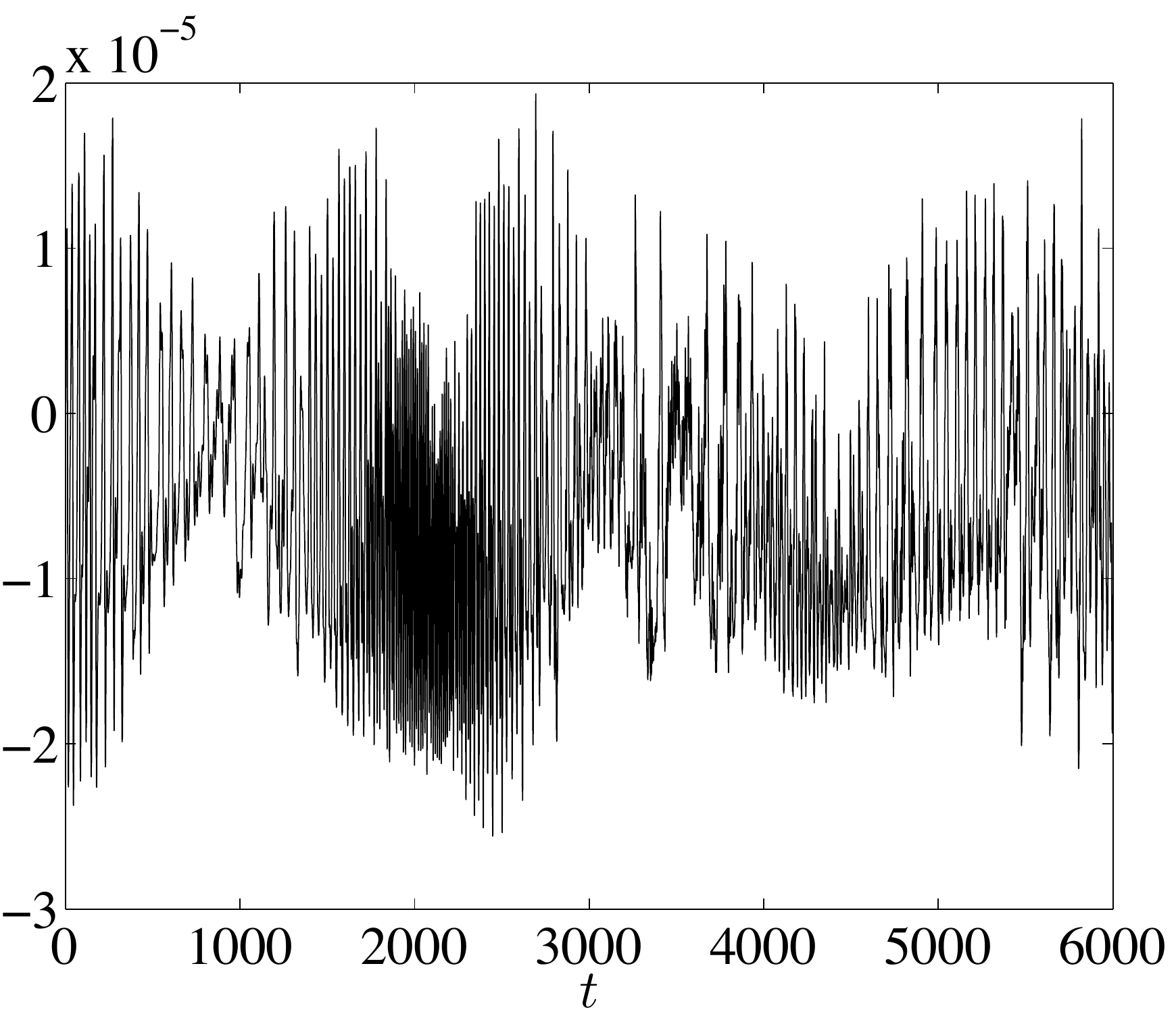}}
}
\centerline{
	\subfigure[Angular velocity of the base spacecraft $\Omega$]{
		\includegraphics[width=0.5\columnwidth]{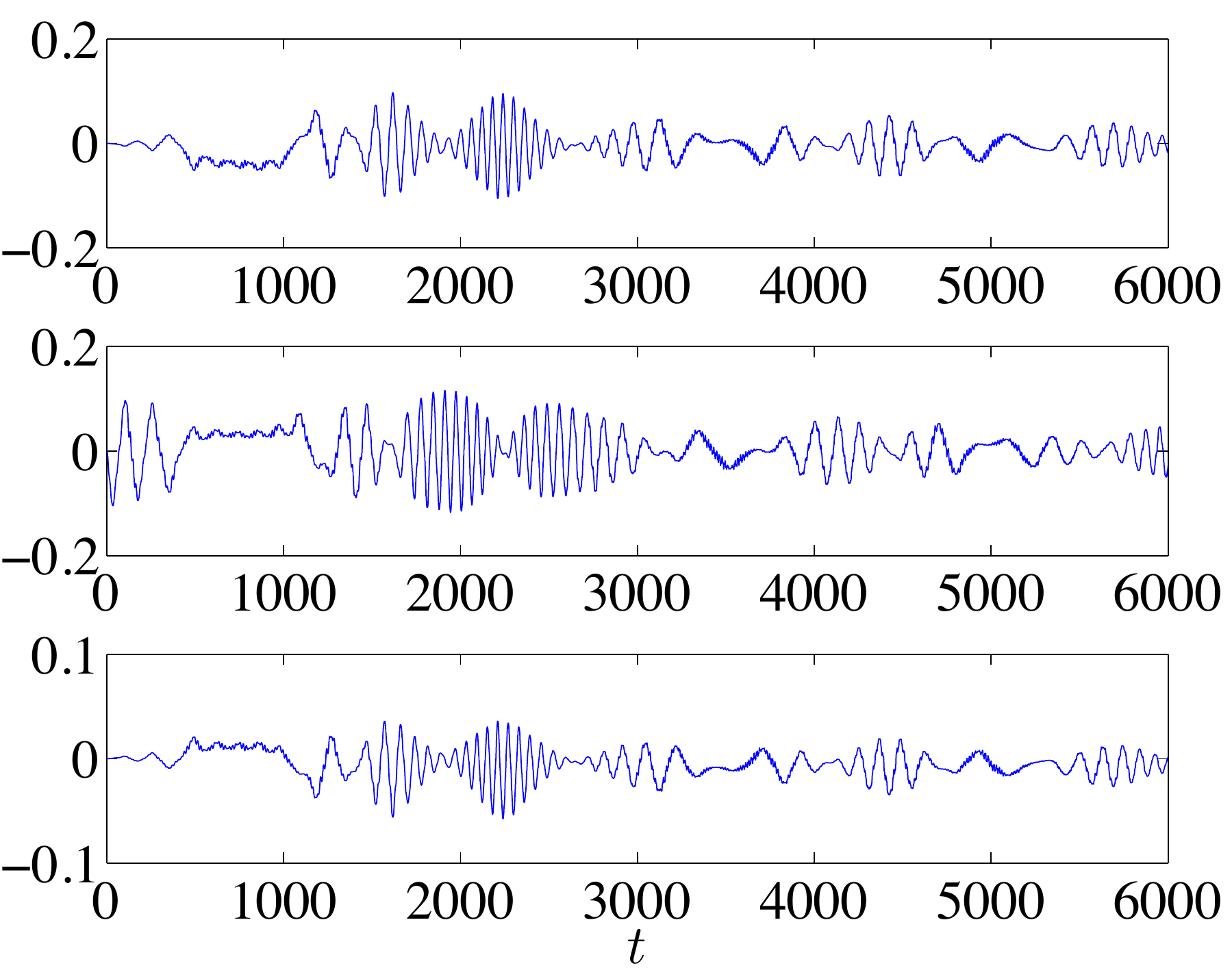}}
	\hspace*{0.00\textwidth}
	\subfigure[Unstretched length of the deployed part of the tether (red), stretched length (blue)]{
		\includegraphics[width=0.52\columnwidth]{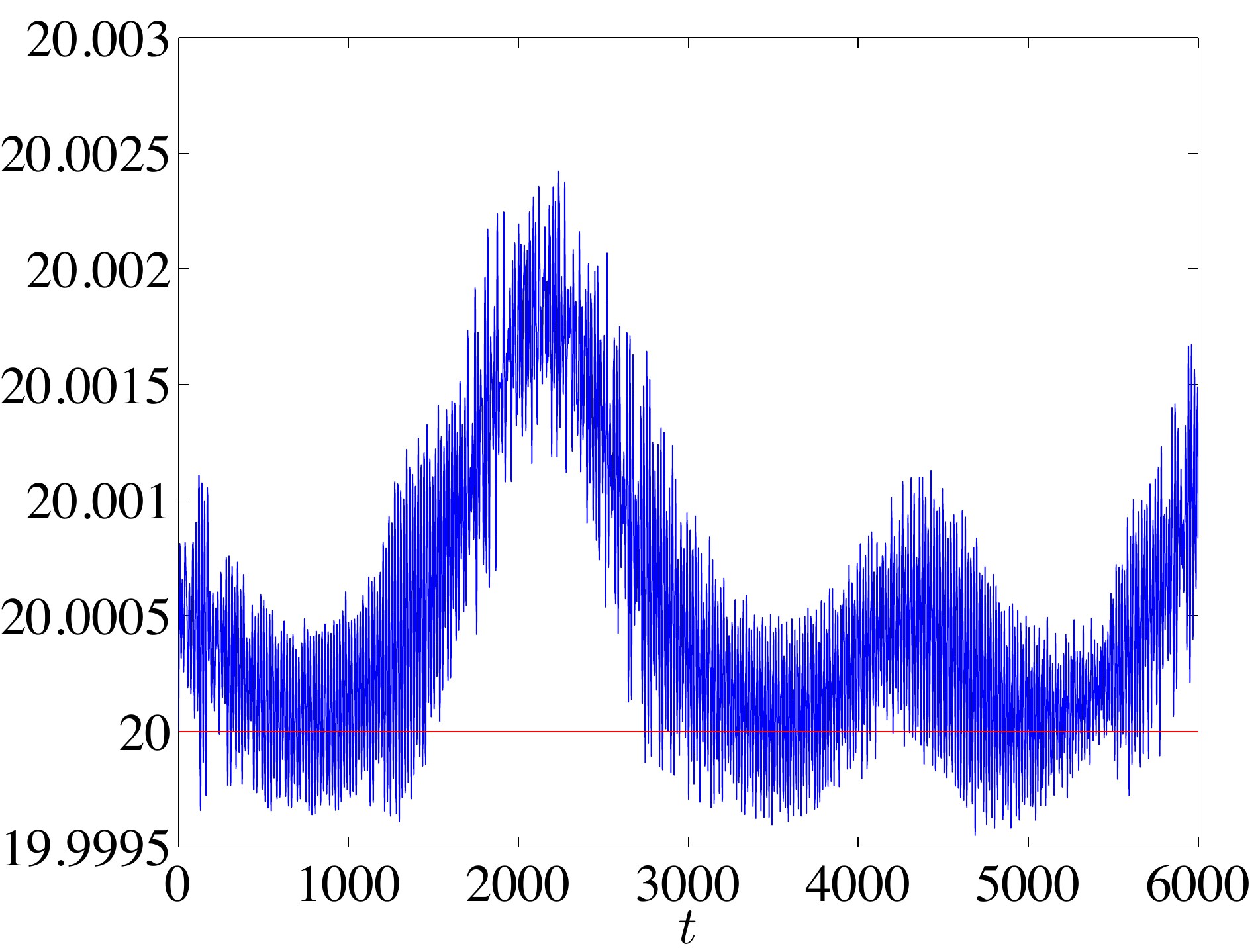}}
}

\caption{Case 1: Circular orbit, Fixed unstretched tether length}
\end{figure}

\begin{figure}
\centerline{
	\subfigure[Snapshots observed at the LVLH frame (km) (The size of spacecraft is increased by a factor of 100 to illustrate attitude dynamics.)]{
		\includegraphics[width=0.90\columnwidth]{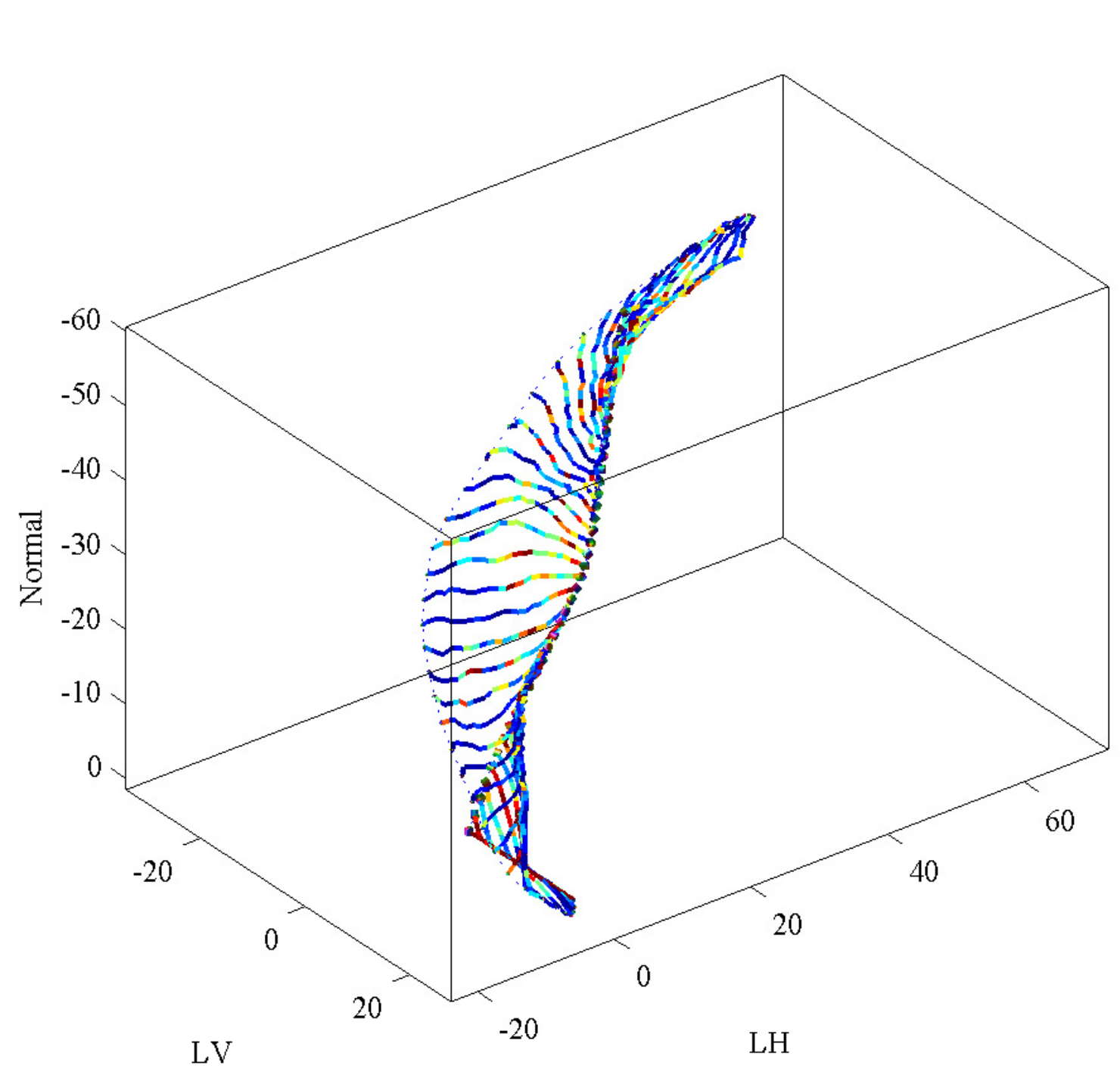}}
}
\centerline{
	\subfigure[$T_{base}+T_{sub}$ (red), $T_{tether}$ (green), $V_{gravity}$ (cyan), $V_{elastic}$ (blue), total energy (black)]{
		\includegraphics[width=0.49\columnwidth]{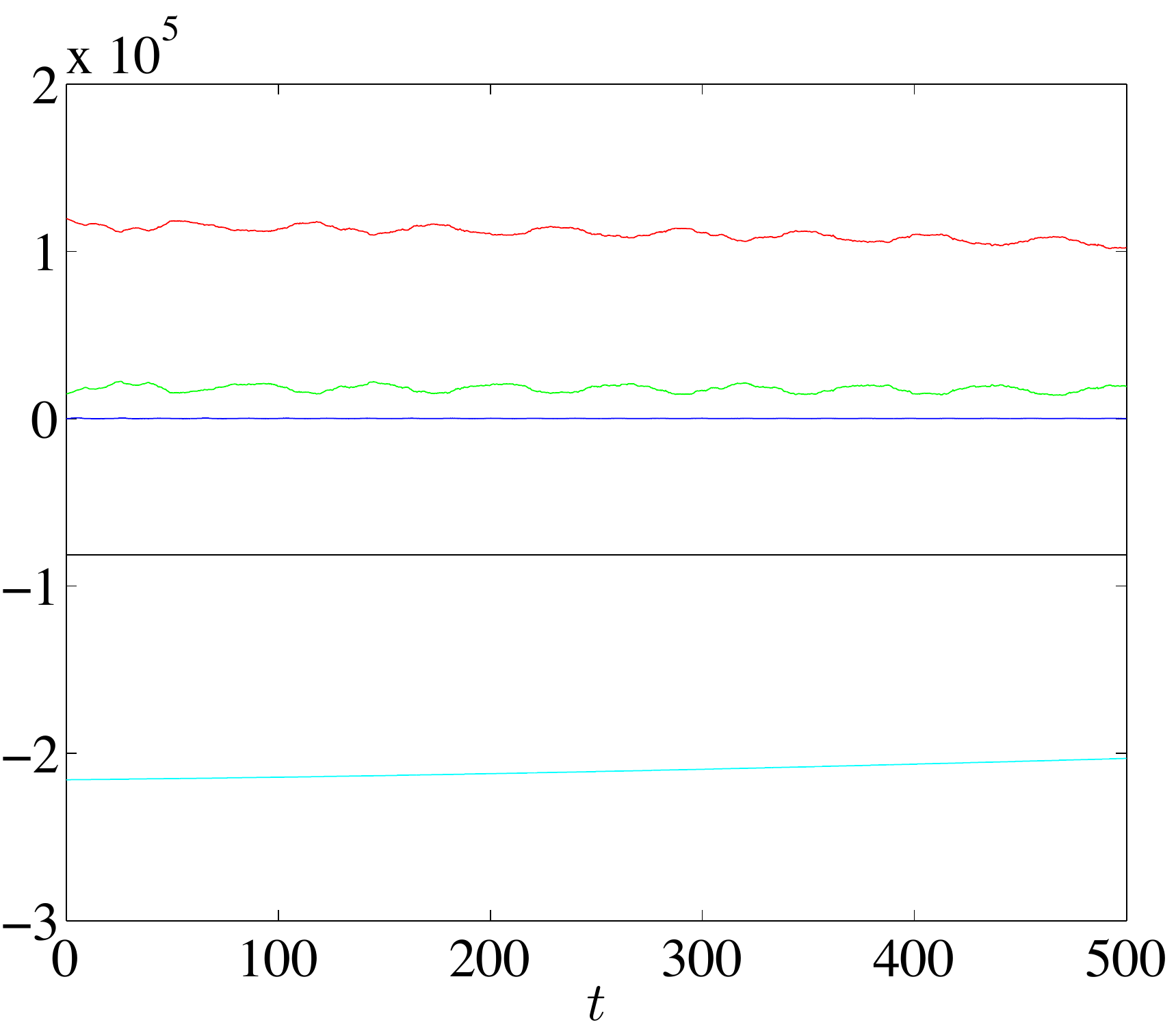}}
	\hspace*{0.00\textwidth}
	\subfigure[Computed total energy deviation $E(t)-E(0)$]{
		\includegraphics[width=0.51\columnwidth]{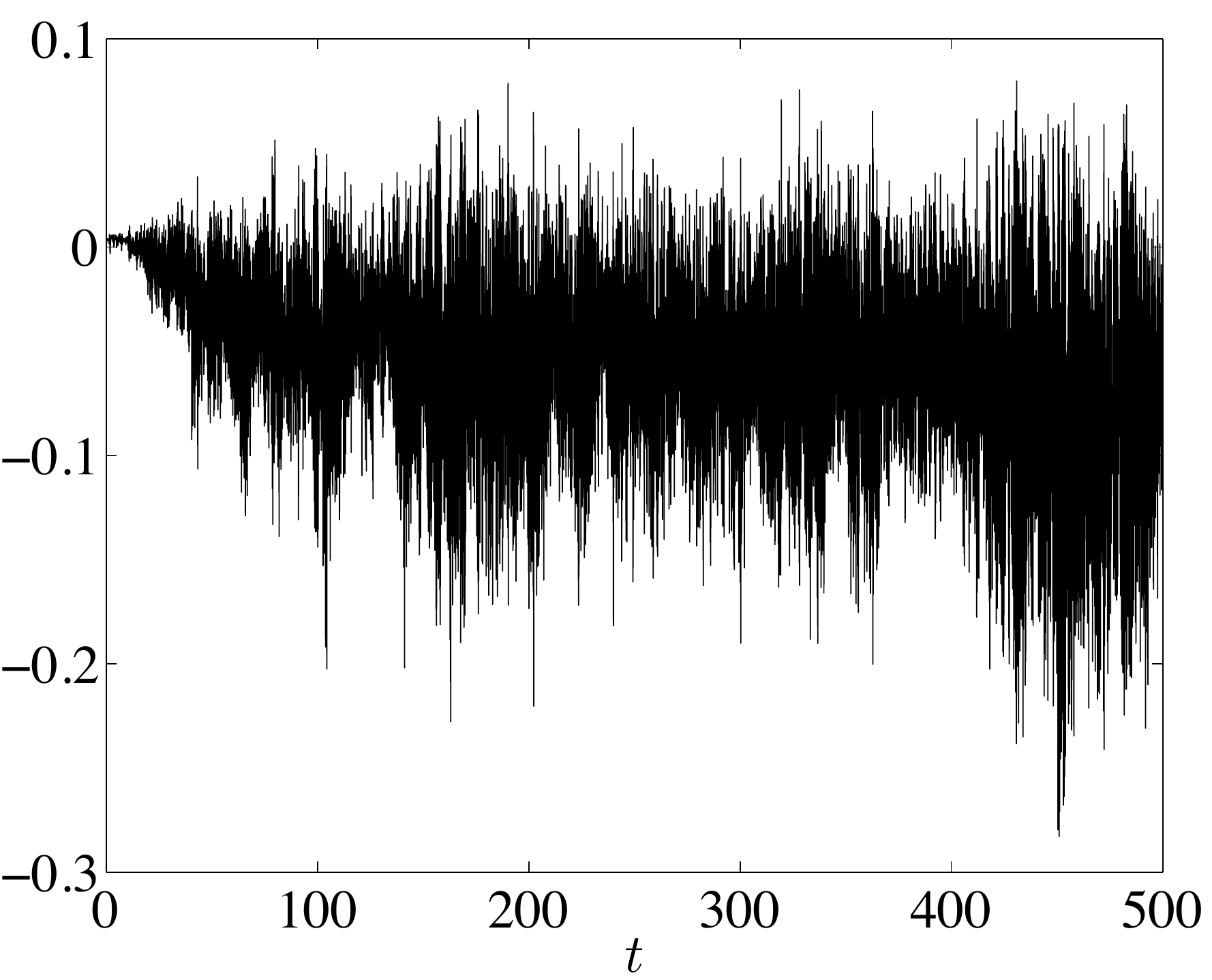}}
}
\centerline{
	\subfigure[Angular velocity of the base spacecraft $\Omega$]{
		\includegraphics[width=0.51\columnwidth]{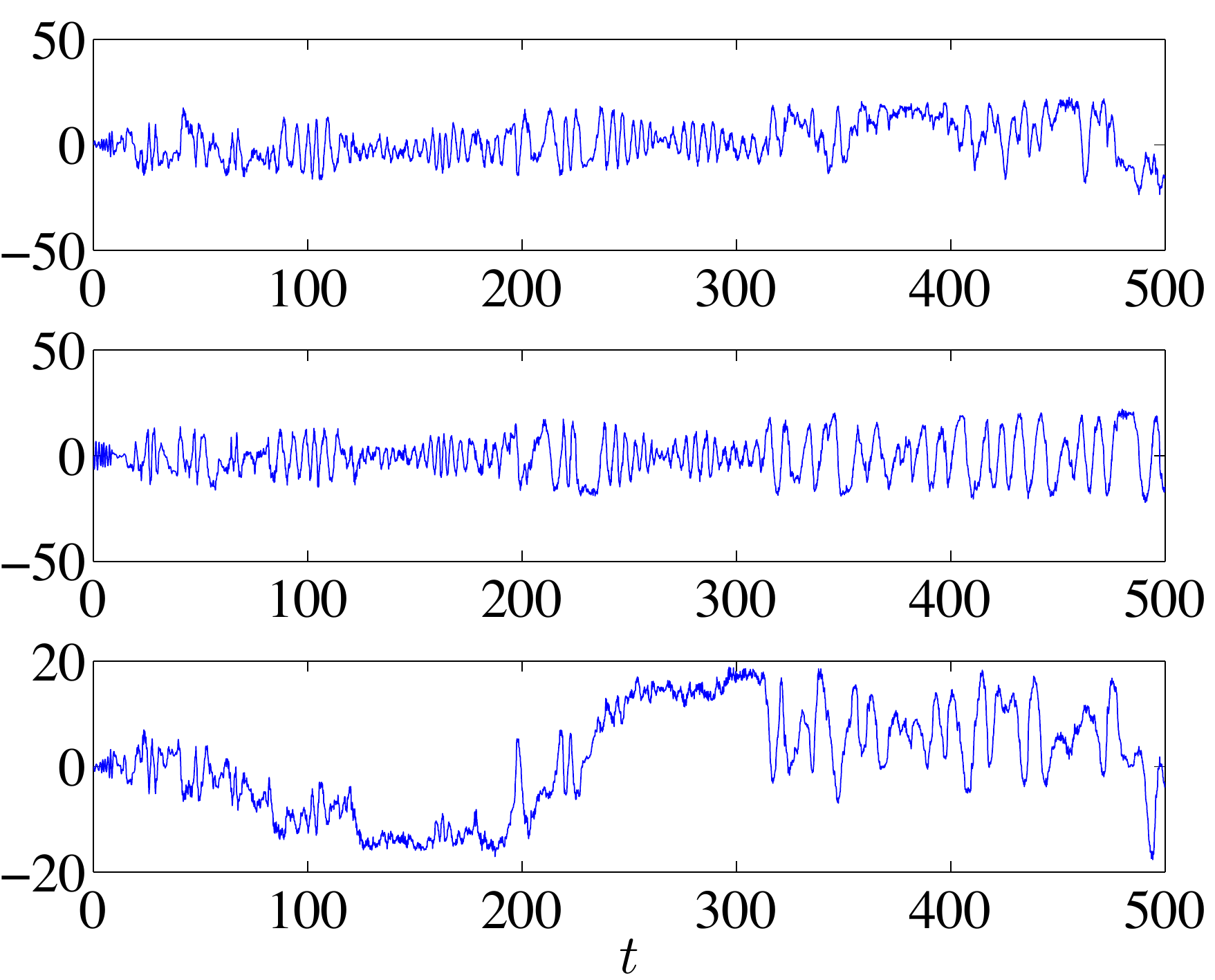}}
	\hspace*{0.00\textwidth}
	\subfigure[Unstretched length of the deployed part of the tether (red), stretched length (blue)]{
		\includegraphics[width=0.49\columnwidth]{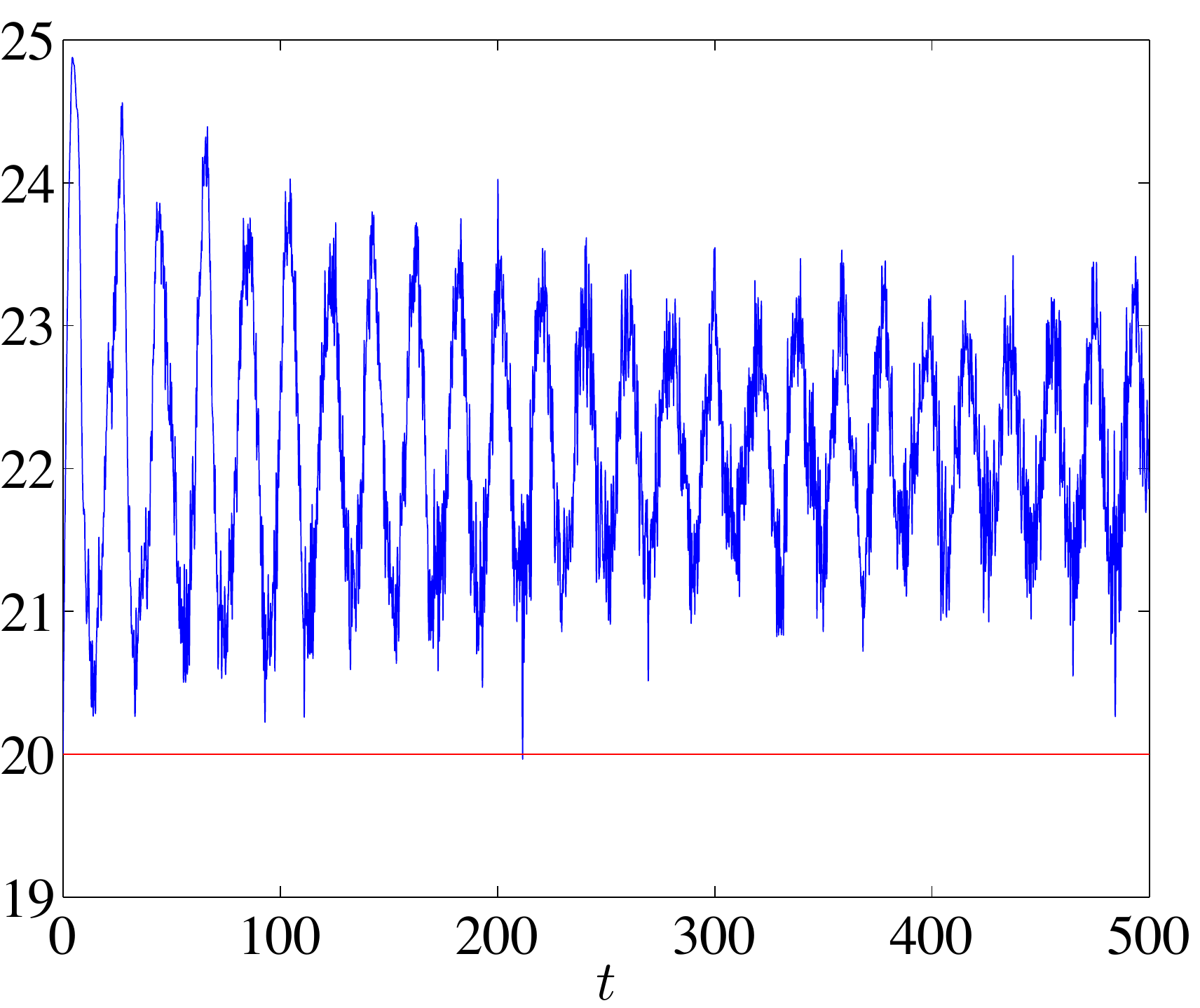}}
}

\caption{Case 3: Perturbed circular orbit, Fixed unstretched tether length}
\end{figure}

\end{document}